\documentclass[12pt,leqno]{article}
\usepackage[british]{babel}
\usepackage[tbtags]{amsmath}
\usepackage{amssymb,amsthm,amsopn,paralist,booktabs}

\theoremstyle{plain}
\newtheorem{theorem}{Theorem}[section]
\newtheorem{proposition}[theorem]{Proposition}
\newtheorem{lemma}[theorem]{Lemma}
\newtheorem{corollary}[theorem]{Corollary}

\theoremstyle{definition}

\theoremstyle{remark}

\newtheorem{remark}[theorem]{Remark}
\newtheorem*{acknowledgements}{Acknowledgements}

\numberwithin{equation}{section}

\newcommand{\Lie}[1]{\operatorname{\textsl{#1}}}
\newcommand{\lie}[1]{\operatorname{\mathfrak{#1}}}

\newcommand{\SO}{\Lie{SO}}
\newcommand{\so}{\lie{so}}

\newcommand{\SU}{\Lie{SU}}
\newcommand{\su}{\lie{su}}
\newcommand{\Un}{\Lie{U}}
\newcommand{\un}{\lie{u}}

\newcommand{\com}{\makebox[7pt]{\raisebox{1.5pt}{\tiny$\circ$}}}

\begin{document}

\title{$\SU(3)$-structures on   hypersurfaces
of  manifolds with $G_2$-structure}

\author{\normalsize\textsc{Francisco Mart\'\i n Cabrera}}

\date{}

\maketitle

\begin{footnotesize}
  \setlength{\parindent}{0pt} \setlength{\parskip}{10pt}


 \textbf{Abstract.} We study $\SU(3)$-structures induced on orientable hypersurfaces
   of seven-dimensional manifolds with $\Lie{G}_2$-structure.
   Taking Gray-Hervella types for both structures
   into account, we relate the type of $\SU(3)$-structure and the
 type of $\Lie{G}_2$-structure with  the shape tensor of the
 hypersurface. Additionally, we show how to compute the intrinsic
 $\Lie{SU}(3)$-torsion and the intrinsic $\Lie{G}_2$-torsion  by means of
 exterior algebra.



\end{footnotesize}

\baselineskip=5mm \pagestyle{myheadings}

\section{Introduction}
The exceptional Lie group $\Lie{G}_2$ is the group of
automorphisms of the Cayley numbers ${\Bbb O}$: the
non-associative normed algebra over ${\Bbb R}$ of dimension eight
having an orthonormal basis $\{ 1 , e_0 , \dots , e_6 \}$ with
multiplication determined by
\begin{gather*}
e_i^2 = -1, \qquad e_i e_j = - e_j e_i, \qquad (i \neq j), \\
e_i e_{i+1} = e_{i+3}, \qquad  e_{i+3} e_{i} = e_{i+1}, \qquad
e_{i+1} e_{i+3} = e_{i}, \qquad (i,j \in {\Bbb Z}_{\;7}).
\end{gather*}
The set $Im{\Bbb O}$ of pure imaginary Cayley numbers is the span
of $\{ e_0, \dots, e_6\}$. The group $\Lie{G}_2$ can be
equivalently defined as the subgroup of $\SO(7)$ acting on
$Im{\Bbb O}$ consisting of those elements which preserve the
three-form given by
\begin{equation} \label{fundphi}
\varphi = \sum_{ i \in {\Bbb Z}_{\;7}} e_i \wedge e_{i+1} \wedge
e_{i+3},
\end{equation}
where we have also denoted by $e_i$ the dual one-form of $e_i$.
Furthermore, it can be shown that $\Lie{G}_2$ acts trivially on
${\Bbb R}$ as a subalgebra of ${\Bbb O}$, so $\Lie{G}_2$ acts
preserving the splitting ${\Bbb O} = {\Bbb R} \oplus
\Lie{Im}\,{\Bbb O}$.

The subgroup of $\Lie{G}_2$ consisting of the automorphisms on
${\Bbb O}$ that leave invariant an imaginary unit, for instance,
$e_0$, is isomorphic to $\SU(3)$. On this fact it is essentially
based the possibility of defining an $\SU(3)$-structure on an
orientable hypersurface of a manifold equipped with a
$\Lie{G}_2$-structure. Thus, Calabi in \cite{Calabi:ac6manifold}
and Gray in \cite{Gray:someexamples} considered orientable
hypersurfaces $M$ of $\Lie{Im}\,{\Bbb O}$, studying the induced
$\SU(3)$-structure on $M$ as an almost Hermitian structure
($\Un(3)$-structure). We recall that manifolds with
$\SU(n)$-structures, i.e., special almost Hermitian manifolds, are
defined as almost Hermitian manifolds $(M,I,\langle \cdot , \cdot
\rangle )$ equipped with a complex volume form $\Psi = \psi_+ + i
\psi_-$. Likewise, a $\Lie{G}_2$-structure on a seven-dimensional
manifold $\overline{M}$ is a reduction of the structure group of
the tangent bundle to $\Lie{G}_2$ and is determined by a global
three-form $\varphi$ which may locally written as in
\eqref{fundphi}.

In the present work,  we consider the structure on $M$ in its full
condition as $\SU(3)$-structure. So, we complete the information
of the particular situations studied by Calabi and Gray with
information coming from the complex volume form $\Psi$. On the
another hand, we extend the results of Calabi and Gray to
orientable hypersurfaces $M$ of manifolds $\overline{M}$ with
$\Lie{G}_2$-structures in general. Thus, in a first instance, we
prove results relating tensors involving the almost complex
structure $I$, the complex volume form $\Psi$, the shape tensor of
$M$ and the three form $\varphi$. Then, fixing the type of
$\Lie{G}_2$-structure on $\overline{M}$, we apply such results in
deriving  geometrical conditions to be satisfied by $M$ to
characterize types of $\SU(3)$-structure.

The paper is organised as follows.  In \S\ref{sec:su3structures}
we present some preliminary material (definitions, results,
notation, etc.) about special almost Hermitian six-manifolds and
give a new alternative way to characterize the diverse types of
such structures (Proposition~\ref{covaseis}). Furthermore, we
explicitly describe the intrinsic $\Lie{SU}(3)$-torsion in terms
of the exterior derivatives of the K{\"a}hler form $\omega$ and
the complex volume form $\Psi$.
 Then in
\S\ref{sec:g2structure}, we recall some definitions, results and
notations about $\Lie{G}_2$-structures. Likewise, we express the
intrinsic $\Lie{G}_2$-torsion  in terms of the exterior
derivatives $d \varphi$ and $d \ast \varphi$.

In  \S\ref{sec:orienthyp} we consider orientable hypersurfaces $M$
of seven-dimensional manifolds $\overline{M}$ with
$\Lie{G}_2$-structure. Thus we start with the description of the
induced $\SU(3)$-structure on $M$. Then we prove the key results
of the present work which are contained in Proposition~\ref{rrB}.
Such results are repeatedly applied in the proves of the
successive theorems. In each one of these theorems, we fix a type
of $\Lie{G}_2$-structure on $\overline{M}$. To avoid a too long
exposition, we only consider types of $\Lie{G}_2$-structure
corresponding to the irreducible modules ${\cal X}_1, \dots ,
{\cal X}_4$ which occurs as summands in the
$\Lie{G}_2$-decomposition of the space ${\cal X}$ of covariant
derivatives of $\varphi$ \cite{Fern-Gray:G2struct}. Finally, in
\S\ref{sec:examples}, we give some examples.

\begin{acknowledgements} This work is supported by a grant from MEC
(Spain), project MTM2004-2644.
\end{acknowledgements}

\section{$\SU(3)$-structures} \label{sec:su3structures}
In this section we give a brief summary of properties and results
relative with special almost Hermitian six-manifolds.  For more
detailed and exhaustive information see
\cite{Chiossi-S:SU3-G2,Cabrera:special}. In the final part of this
section, we present a new alternative way to characterize the
diverse types of $\SU(3)$-structure (Proposition~\ref{covaseis})
and express the intrinsic $\Lie{SU}(3)$-torsion by means of
exterior algebra.

An {\it almost Hermitian} manifold is a $2n$-dimensional manifold
$M$ with a $\Un(n)$-structure. This means that $M$ is equipped
with a Riemannian metric $\langle \cdot , \cdot \rangle$ and an
orthogonal almost complex structure $I$. Each fibre $T_m M$ of the
tangent bundle  can be consider as complex vector space by
defining $i x = Ix$.  We will write $T_m M_{\Bbb C}$ when we are
regarding $T_m M$ as such a space.

We define a Hermitian scalar product $\langle \cdot , \cdot
\rangle_{\Bbb C} = \langle \cdot , \cdot  \rangle + i
\omega(\cdot,\cdot)$, where $\omega$ is the K{\"a}hler form given
by $\omega (x, y) = \langle x, Iy \rangle$. The real tangent
bundle $TM$ is identified with the cotangent bundle $T^*M$ by the
map $x \to \langle \cdot , x \rangle=x$. Analogously, the
conjugate complex vector space $\overline{T_m M}_{\Bbb C}$ is
identified with the dual complex space $T^{\ast}_m M_{\Bbb C}$ by
the map $x \to \langle \cdot , x \rangle_{\Bbb C} = x_{\Bbb C}$.
It follows immediately that $x_{\Bbb C} = x + i Ix$.

If we consider the spaces $\Lambda^p T^*_m M_{\Bbb C}$ of
skew-symmetric complex forms, one can check  $x_{\Bbb C} \wedge
y_{\Bbb C} = (x + i Ix) \wedge (y + i Iy)$. There are natural
extensions of scalar products to $\Lambda^p T^*_m M$ and
$\Lambda^p T^*_m M_{\Bbb C}$, defined respectively by
\begin{eqnarray*}
\displaystyle \langle a,b\rangle & = &  \frac{1}{p!} \sum_{i_1,
\dots , i_p=1}^{2n}
a(e_{i_1}, \dots, e_{i_p}) b(e_{i_1},\dots, e_{i_p}), \\
\displaystyle  \langle a_{\Bbb C} ,b_{\Bbb C} \rangle_{\Bbb C} & =
& \frac{1}{p!} \sum_{i_1, \dots , i_p=1}^{n} a_{\Bbb C} (u_{i_1},
\dots, u_{i_p}) \overline{ b_{\Bbb C}(u_{i_1},\dots, u_{i_p})},
\end{eqnarray*}
where $e_1, \dots , e_{2n}$ is an orthonormal basis for real
vectors and $u_1, \dots, u_n$ is a unitary basis for complex
vectors.

 The following conventions will be
used in this paper.  If $b$ is a $(0,s)$-tensor, we write
\begin{equation*} \label{ecuacionesb}
  \begin{array}{l}
    I_{(i)}b(X_1, \dots, X_i, \dots , X_s) = - b(X_1, \dots , IX_i, \dots ,
    X_s),\\[2mm]
    I b(X_1,\dots,X_s) = (-1)^sb(IX_1,\dots,IX_s),
  \end{array}
\end{equation*}

A {\it special almost Hermitian} manifold is a $2n$-dimensional
manifold $M$ with an $\SU(n)$-structure. This means that $(M,
\langle \cdot, \cdot \rangle, I)$ is an almost Hermitian manifold
equipped with a complex volume form $\Psi = \psi_+ + i \psi_-$
such that $\langle  \Psi , \Psi \rangle_{\Bbb C} = 1$. Note that
$I_{(i)} \psi_+ = \psi_-$. \vspace{2mm}

In the following, we will only consider special almost Hermitian
six-manifold.
 If  $e_1,e_2,e_3$ is a unitary basis for complex vectors such that \linebreak $  \Psi
(e_1, e_2, e_3)=1$,  i.e., $ \psi_+ (e_1, e_2 , e_3) = 1$ and
$\psi_- (e_1 , e_2 , e_3)=0$,  then \linebreak $e_1$, $e_2$,
$e_3$, $Ie_1$, $Ie_2$, $Ie_3$ is an orthonormal basis for real
vectors {\it adapted} to the $\SU(3)$-structure. The three forms
$\psi_+$ and $\psi_-$ can be  respectively expressed by
\begin{gather*}
\psi_+ = e_1 \wedge e_2 \wedge e_3 - Ie_1 \wedge Ie_2 \wedge e_3
-Ie_1 \wedge e_2 \wedge Ie_3 -e_1 \wedge Ie_2 \wedge Ie_3, \\
\psi_- = - Ie_1 \wedge Ie_2 \wedge Ie_3 + Ie_1 \wedge e_2 \wedge
e_3 + e_1 \wedge Ie_2 \wedge e_3 + e_1 \wedge e_2 \wedge Ie_3.
\end{gather*}
Moreover, it is straightforward to check $
 \omega^3 = 6 \, e_1 \wedge e_2 \wedge e_3 \wedge Ie_1 \wedge Ie_2 \wedge Ie_3,
$
 where $\omega^3 = \omega \wedge \omega \wedge  \omega$. If we fix the form $Vol$
 such that $6 Vol = \omega^3$ as real volume form, it
follows next lemma.
\begin{lemma}[\cite{Cabrera:special}] \label{volumenes}
Let $M$ be a special almost Hermitian $6$-manifold, then
\begin{enumerate}
\item[{\rm (i)}] $\psi_+ \wedge \omega = \psi_- \wedge \omega = 0
$;
 \item[{\rm (ii)}] $ \psi_+ \wedge \psi_- = -   4 \, Vol$ and $
\psi_+ \wedge \psi_+ = \psi_- \wedge \psi_- = 0$;
 \item[{\rm (iii)}] for  $1 \leq i < j \leq 3$, $I_{(i)} I_{(j)} \psi_+ = - \psi_+$
and $I_{(i)} I_{(j)} \psi_- = - \psi_-$; and
 \item[{\rm (v)}] $ x \wedge \psi_+ =  Ix \wedge \psi_-= - (Ix \lrcorner \psi_+) \wedge
\omega$ and  $x \lrcorner \psi_+ = Ix \lrcorner \psi_-$, for all
vector $x$, where $\lrcorner$ denotes the interior product.
\end{enumerate}
\end{lemma}

 Furthermore,
let us note that there are two Hodge star operators, defined on
$M$, respectively associated with the volume forms $Vol$ and
$\Psi$. Relative to the real Hodge star operator $\ast$, we have
the equations
\begin{gather} \label{estrella1}
 \ast \left( \ast ( \mu \wedge \psi_+ ) \wedge \psi_+ \right) =
\ast \left( \ast ( \mu \wedge \psi_- ) \wedge \psi_- \right)  = -2
\mu, \\
 \ast \left( \ast ( \mu \wedge \psi_- ) \wedge
\psi_+ \right) = - \ast \left( \ast ( \mu \wedge \psi_+ ) \wedge
\psi_- \right) =  2
I\mu. \label{estrella2}
\end{gather}

We are dealing with $\Lie{G}$-structures where $\Lie{G}$ is a
subgroup of the linear group $\Lie{GL}(m , {\Bbb R})$. If $M$
possesses a $\Lie{G}$-structure, then there always exists a
$\Lie{G}$-connection defined on $M$. Moreover, if $(M^m ,\langle
\cdot , \cdot \rangle)$ is an orientable $m$-dimensional
Riemannian manifold and  $\Lie{G}$ a closed and connected subgroup
of $\SO(m)$, then there exists a unique metric $G$-connection
$\widetilde{\nabla}$ such that $\xi_x = \widetilde{\nabla}_x -
\nabla_x$ takes its values in $\lie{g}^{\perp}$, where
$\lie{g}^{\perp}$ denotes the orthogonal complement in $\so(m)$ of
the Lie algebra $\lie{g}$ of $\Lie{G}$ and $\nabla$ denotes the
Levi-Civita connection
\cite{Salamon:holonomy,CleytonSwann:torsion}. The tensor $\xi$ is
the {\it intrinsic torsion}  of the $\Lie{G}$-structure and
$\widetilde{\nabla}$ is called the {\it minimal
$\Lie{G}$-connection}.

 For $\Un(3)$-structures, the minimal $\Un(3)$-connection is given by
$ \widetilde{\nabla} = \nabla + \xi$, with
\begin{equation} \label{torsion:xi}
\xi_X Y = - \frac{1}{2} I \left( \nabla_X I \right) Y.
\end{equation}
see \cite{Falcitelli-FS:aH}. Since $\Un(3)$ stabilises the
K{\"a}hler form $\omega$, it follows that $\widetilde{\nabla}
\omega = 0$. Then $\nabla \omega = - \xi \omega \in T^* M \otimes
\un(3)^{\perp}$. Thus, one can identify the $\Un(3)$-components of
$\xi$ with the $\Un(3)$-components of $\nabla \omega$.

For $\SU(3)$-structures, we have the decomposition $\lie{so}(6)=
\lie{su}(3) + {\Bbb R} + \un(3)^{\perp}$, i.e.,
$\lie{su}(3)^{\perp} = {\Bbb R} + \un(3)^{\perp}$. Therefore, the
intrinsic $\SU(3)$-torsion $\eta + \xi$ is such that $\eta \in T^*
M \otimes {\Bbb R} \cong T^* M$ and $\xi$ is still determined by
Equation \eqref{torsion:xi}. The tensors $\omega$, $\psi_+$ and
$\psi_-$ are stabilised by the $\SU(3)$-action, and
$\widetilde{\nabla}^{\SU(3)} \omega = 0$,
$\widetilde{\nabla}^{\SU(3)} \psi_+ = 0$ and $
\widetilde{\nabla}^{\SU(3)} \psi_- = 0$, where
$\widetilde{\nabla}^{\SU(3)} = \nabla + \eta + \xi$ is the minimal
$\SU(3)$-connection. Since $\widetilde{\nabla}^{\SU(3)}$ is metric
and $\eta \in T^* M \otimes {\Bbb R}$, we have $ \langle Y ,
\eta_X Z \rangle = (I \eta)(X) \omega(Y,Z)$, where $\eta$ on the
right side is  a one-form. Hence
 \begin{equation} \label{torsion:eta} \eta_X Y = I
\eta (X) IY.
\end{equation}

One can check  $\eta \omega = 0$, then  from
$\widetilde{\nabla}^{\SU(3)} \omega = 0$ it is obtained
$$
 \nabla \omega = - \xi \omega \in T^* M
\otimes \un(3)^{\perp} = {\cal W}_1 + {\cal W}_2 + {\cal W}_3 +
{\cal W}_4,
$$
where the summands ${\cal W}_i$ are the irreducible
$\Un(3)$-modules given by Gray and Hervella \cite{Gray-H:16} and
$+$ denotes direct sum. The spaces ${\cal W}_3$ and ${\cal W}_4$
are also irreducible as $\SU(3)$-modules. However, ${\cal W}_1$
and ${\cal W}_2$ admit the  decompositions $ {\cal W}_i = {\cal
W}_i^+ + {\cal W}_i^-$, $i=1,2$, into irreducible
$\SU(3)$-components, where ${\cal W}_i^+$ (${\cal W}_i^-$)
includes those  $a \in {\cal W}_i \subseteq T^*M \otimes \Lambda^2
T^* M$ such that the bilinear form $r(a)$, defined by $2 r(a)(x,y)
= \langle  x \lrcorner a , y \lrcorner \psi_+ \rangle$, is symmetric
 (skew-symmetric). 

On the other hand, since $\widetilde{\nabla}^{\SU(3)} \psi_+ = 0$
and $\widetilde{\nabla}^{\SU(3)} \psi_- = 0$, we have  $ \nabla
\psi_+ = - \eta \psi_+ - \xi \psi_+$ and $ \nabla \psi_- = - \eta
\psi_- - \xi \psi_-$. Therefore, from Equations \eqref{torsion:xi}
and \eqref{torsion:eta} we obtain the following expressions
\begin{equation*} \label{torsiones}
\begin{array}{lll}
\quad - \eta_X \psi_+   =  - 3 I\eta(X) \psi_-, & \quad &  -\xi_X
\psi_+  =
\displaystyle \frac{1}{2}  (e_i \lrcorner \nabla_X \omega) \wedge
(e_i \lrcorner \psi_-), \\[2mm]
\quad - \eta_X \psi_-   =    3 I\eta(X) \psi_+,  & &
 -\xi_X \psi_-  =
 - \displaystyle \frac{1}{2}  (e_i \lrcorner \nabla_X \omega)
\wedge (e_i \lrcorner \psi_+),
\end{array}
\end{equation*}
where the summation convention is used. We use this convention
along the present paper. 

 It is obvious that $- \eta \psi_+ \in {\cal W}^-_{5}
 = T^* M \otimes \psi_-$  and
 $- \eta \psi_- \in {\cal W}^+_{5} =  T^* M \otimes \psi_+$.
The tensors $-\xi \psi_+$ and $-\xi \psi_-$ are described in the
following proposition, where we need to consider the two
$\SU(3)$-maps
\begin{equation} \label{Xiplusminus}
 \Xi_+ , \Xi_- \, : \, T^* M \otimes {\frak
u}(3)^{\perp} \to T^* M \otimes \Lambda^3 T^* M
\end{equation}
respectively defined by $\nabla_{\cdot} \omega \to 1/2 \, (e_i
\lrcorner \nabla_{\cdot} \omega) \wedge (e_i \lrcorner \psi_-)$
 and
 $ \nabla_{\cdot} \omega \to \linebreak  - 1/2 \,
 (e_i \lrcorner \nabla_{\cdot} \omega) \wedge (e_i \lrcorner \psi_+)$.
\begin{proposition}[\cite{Cabrera:special}] \label{Ximasmenos}
The $SU(3)$-maps $\Xi_+$ and $\Xi_-$ are injective and
\begin{equation*}
\Xi_{+} \left( T^* M \otimes \un(3)^{\perp} \right) = \Xi_{-}
\left( T^* M \otimes \un(3)^{\perp} \right) = T^* M \otimes T^* M
\wedge \omega.
\end{equation*}
\end{proposition}

 Last result  and above considerations give rise to the  following
theorem which describes the tensors $\nabla \psi_+$ and $\nabla
\psi_-$.
\begin{theorem}[\cite{Cabrera:special}]
 \label{igual3} Let $M$ be a special almost Hermitian
 $6$-manifold
 with K{\"a}hler form $\omega$ and complex volume form $\Psi =
\psi_+ + i \psi_-$. Then
$$
\begin{array}{l}
\nabla \psi_+ \in {\cal W}^{\Xi ; a}_1 + {\cal W}^{\Xi ; b}_1 +
{\cal W}^{\Xi ; a}_2 +  {\cal W}^{\Xi ; b}_2 + {\cal
W}^{\Xi}_3 + {\cal W}^{\Xi}_4 + {\cal W}^-_5, \\[2mm]
 \nabla \psi_-
\in {\cal W}^{\Xi ; a }_1 + {\cal W}^{\Xi ; b}_1  + {\cal W}^{\Xi
; a}_2 + {\cal W}^{\Xi ; b }_2 + {\cal W}^{\Xi}_3 + {\cal
W}^{\Xi}_4 + {\cal W}^+_5,
\end{array}
$$
where ${\cal W}^{\Xi;a}_i = \Xi_{+} ({\cal W}^+_i) = \Xi_{-}
({\cal W}^-_i)$, ${\cal W}^{\Xi;b}_i = \Xi_{+} ({\cal W}^-_i) =
\Xi_{-} ({\cal W}^+_i)$, $i=1,2$; ${\cal W}^{\Xi}_j = \Xi_{+}
({\cal W}_j) = \Xi_{-} ({\cal W}_j)$, $j=3,4$;
 ${\cal W}_{5}^+ = T^*M \otimes \psi_{+}$ and
${\cal W}_{5}^- = T^*M \otimes \psi_{-}$.
\end{theorem}
Further details of the $\SU(3)$-components of  $\nabla \psi_+$ and
$\nabla \psi_-$ can be found in \cite{Cabrera:special}. Moreover,
if we consider the alternation maps $\widetilde{\it a}_{\pm} \, :
\, {\cal W}^{\Xi} + {\cal W}_5^{\mp} \to \Lambda^{4} T^* M$, we
get the following consequences of Theorem \ref{igual3}.
\begin{corollary}[\cite{Cabrera:special}] \label{corigual3}
For  $\SU(3)$-structures, the exterior derivatives of $\psi_+$ and
$\psi_-$ are such that
$$
d \psi_+, d \psi_- \in \Lambda^{4} T^* M = {\cal W}_1^{a} + {\cal
W}_2^{a} + {\cal W}_{4,5}^{a},
$$
where  $\widetilde{a}_{\pm} ({\cal W}_1^{\Xi;b}) = {\cal
W}_1^{a}$, $\widetilde{a}_{\pm}({\cal W}_2^{\Xi;b}) ={\cal
W}_2^{a}$ and $\widetilde{a}_{\pm}({\cal W}_4^{\Xi}) =
\widetilde{a}_{\pm}({\cal W}_5^{\mp}) = {\cal W}_{4,5}^{a}$.
Moreover, $\Lie{Ker}(\widetilde{a}_{\pm}) = {\cal W}_1^{\Xi;a} +
{\cal W}_2^{\Xi;a} + {\cal W}^{\Xi}_3 + {\cal A}_{\pm}$, where
$T^* M \cong {\cal A}_{\pm} \subseteq {\cal W}_4^{\Xi} + {\cal
W}_5^{\mp}$, and the modules ${\cal W}_i^{a}$ are described by
${\cal W}_1^{a} = {\Bbb R}  \omega \wedge \omega$, ${\cal W}_2^{a}
=  {\frak s}{\frak u}(3) \wedge \omega$ and $ {\cal W}_{4,5}^{a} =
T^*M \wedge \psi_+ = T^* M \wedge \psi_- = \un(3)^{\perp} \wedge
\omega$.
\end{corollary}
Now, if we  compute the images by the maps $\widetilde{a}_{\pm}
\com \Xi_{\pm}$ of the ${\cal W}_4$-part of $\nabla \omega$, we
obtain the ${\cal W}_{4,5}^a$-parts of $d \psi_+$ and $d \psi_-$,
i.e.,
\begin{equation} \label{dcuatrocinco}
(d \psi_{\pm})_{4,5} = - \left( 3 \eta + \frac12 I d^* \omega
\right) \wedge \psi_{\pm},
\end{equation}
where $d^* \omega$ means the coderivative of $\omega$. Then,
making use of Equations \eqref{estrella1} and \eqref{estrella2} in
Equation \eqref{dcuatrocinco}, one can   explicitly describe the
one-form $\eta$. This will complete the definition of
 the $\SU(3)$-connection $\widetilde{\nabla}^{\SU(3)}$.
\begin{theorem}[\cite{Cabrera:special}] \label{torsion3w5}
For an $\SU(3)$-structure,
 the ${\cal W}_5$-part $\eta$ of the torsion can be
identified with $- \eta \psi_+ = - 3 I \eta \otimes \psi_-$ or $-
\eta \psi_- =  3 I \eta \otimes \psi_+$,
 where $\eta$ is a one-form such that
$$
\ast \left( \ast  d \psi_{\pm}  \wedge  \psi_{\pm} \right) = 6
\eta +  I d^* \omega = - I \ast \left( \ast  d \psi_+ \wedge
\psi_- \right) =    I \ast \left( \ast  d \psi_- \wedge  \psi_+
\right).
$$
\end{theorem}

Finally, we are going to show an alternative way to describe the
summand $\xi$ of the intrinsic torsion  of an $\SU(3)$-structure.
In fact, if $e_1, e_2 , e_3, e_4= Ie_1, e_5=Ie_2, e_6= Ie_3$ is an
adapted basis to an $\SU(3)$-structure, we have
\begin{equation} \label{covaseis}
\nabla \omega = \sum_{i,j=1}^6 a_{ij} e_i \otimes e_j \lrcorner
\psi_+.
\end{equation}
Now, we consider the $\SU(3)$-map $r : T^* M \otimes
\un(3)^{\perp} \to \otimes^2 T^*M$ defined by
\begin{equation} \label{erresu}
 r(\alpha)(x,y) = \frac{1}{2} \langle x
\lrcorner \alpha , y \lrcorner \psi_+   \rangle.
\end{equation}
It is straightforward to check that, for $\alpha$ given by
Equation~\eqref{covaseis}, \linebreak $r(\alpha) = \sum_{i,j=1}^6
a_{ij} e_i \otimes e_j$ and the following results follow.
\begin{proposition} \label{suseis}
 If $r$ is the map defined by
\eqref{erresu}, then
\begin{enumerate}
\item[{\rm (i)}]
 $r$ is an $\SU(3)$-isomorphism;
 \item[{\rm (ii)}]
 $ r({\cal W}_1^+) = {\Bbb R}\langle \cdot , \cdot \rangle$,
$ r({\cal W}_1^-) = {\Bbb R} \omega$;
 \item[{\rm (iii)}]
 $ r({\cal W}_2^+) = \{ b \in
\otimes^2 T^*M \, | \, b \mbox{ is trace free symmetric and } Ib=b
\}$; \item[{\rm (iv)}]
 $ r({\cal W}_2^-) = \{ b \in \otimes^2 T^*M \, | \, b
\mbox{ is skew-symmetric, }  Ib=b \mbox{ and } \langle \omega , b
\rangle = 0 \,\}$;
 \item[{\rm (v)}] $ r({\cal
W}_3) = \{ b \in \otimes^2 T^*M \, | \, b \mbox{ is symmetric and
} Ib=-b \}$; and
 \item[{\rm (vi)}] $ r({\cal W}_4) = \{ b \in
\otimes^2 T^*M \, | \, b \mbox{ is skew-symmetric and } Ib=-b \}$.
\end{enumerate}
\end{proposition}
Finally, we point out that the coderivative of $\omega$ is given
by
\begin{equation} \label{codeseis}
d^* \omega = \sum_{i=1}^{6} \sum_{\{j,k \, | \,
\psi_+(e_i,e_j,e_k)=1\}} \left(r(\omega)(e_j,e_k)-
r(\omega)(e_k,e_j)\right) e_i,
\end{equation}
where we denote $r(\omega) = r(\nabla \omega)$.

In \cite{Cabrera:special} it is proved that the intrinsic
$\Lie{SU}(3)$-torsion $\eta + \xi$ can be computed from  $d
\omega$, $d \psi_+$ and $d \psi_-$. An explicit description of
this fact is given in next Theorem, where we will write $(d
\psi_{\pm})_{\xi} = d \psi_{\pm} + 3 \eta \wedge \psi_{\pm}$, and
 $(X\wedge Y)\lrcorner d \psi_{\pm}= d\psi_{\pm}(X,Y,\cdot,\cdot
)$. We recall that $\eta$ is given in Theorem \ref{torsion3w5}, it
only remains to describe $\xi$.
\begin{theorem}
 The $\lie{u}(3)$-part $\xi$ of the
torsion of the minimal $\Lie{SU}(3)$-connection
$\widetilde{\nabla}^{\Lie{SU}(3)} = \nabla + \eta + \xi$ is given
by
\begin{equation} \label{exptorsion:xi}
 \xi_{X} Y = - \frac{1}{2} r(\omega)(X, e_i) \psi_{+}
(e_i,e_j,Y) Ie_j,
\end{equation}
for all vectors $X,Y$. Furthermore, the $(0,2)$-tensor field
$r(\omega)$ can be alternatively given by
\begin{gather} \label{erreextalg}
2 r(\omega)(X,Y) = \langle X \lrcorner d \omega , Y \lrcorner
\psi_+ \rangle - \langle (X \wedge Y) \lrcorner (d
\psi_+)_{\xi},\omega \rangle \\
\vspace{6cm} \quad \qquad  + \langle (IX \wedge Y) \lrcorner (d
\psi_-)_{\xi},\omega \rangle. \nonumber
\end{gather}
\end{theorem}
\begin{proof} From Equation \eqref{covaseis} it follows
\begin{equation} \label{dij}
d \omega = a_{ij} e_i \wedge (e_j \lrcorner \psi_+).
\end{equation}
Furthermore, computing $\Xi_{+} (\nabla_{\cdot} \omega)$ and
$\Xi_{-}(\nabla)_{\cdot}\omega)$, defined by Equation
\eqref{Xiplusminus}, and then alternating the obtained result, we
will get
\begin{gather} \label{psimasij}
(d \psi_+)_{\xi} = a_{ij} e_i \wedge e_{j} \wedge \omega, \\
(d \psi_-)_{\xi} = - a_{ij} e_i \wedge Ie_{j} \wedge \omega.
\label{psimenosij}
\end{gather}
From Equations \eqref{dij}, \eqref{psimasij} and
\eqref{psimenosij}, one can check Equation \eqref{erreextalg}.
Finally, from Equation \eqref{torsion:xi} we can  obtain
$$
\xi_X Y = - \frac{1}{2} (\nabla_X \omega)(e_j, Y) Ie_j.
$$
Now using Equation \eqref{covaseis}, we will obtain Equation
\eqref{exptorsion:xi}.
\end{proof}

\section{$\Lie{G}_2$-structures}\label{sec:g2structure}

   We recall briefly some facts about $\Lie{G}_2$-structures and
derive some results we will need in subsequent sections.

A $\Lie{G}_2$-structure on a Riemannian seven-manifold
$(\overline{M},\langle \cdot , \cdot \rangle)$ is by definition a
reduction of the structure group of the tangent bundle to
$\Lie{G}_2$. This is equivalent to the existence of a global
three-form $\varphi$ which may be locally written as in
\eqref{fundphi}. For all $m \in \overline{M}$, the tangent space
$T\,\overline{M}$ is then associated to the representation
$\Lie{Im}{\Bbb O}$ of $\Lie{G}_2$.   Since $\Lie{G}_2$ preserves
$\langle \cdot , \cdot \rangle$, one can define a two-fold vector
cross product $P$ given by $\langle P(x,y),z \rangle = \varphi
(x,y,z)$.  A local orthonormal frame  $\{ e_0, \dots , e_6 \}$ for
vectors is  a {\it Cayley frame}, if $P(e_i,e_{i+1}) = e_{i+3}$
for all $i \in {\Bbb Z}_{\;7}$
\cite{Fern-Gray:G2struct,Cabrera:G2struct}. The four-form $\ast
\varphi$ defined by $\ast \varphi (x,y,z,u) = \lie{S}_{yzu}
\langle P(x,y) , P(z,u) \rangle$, $\lie{S}$ denotes ciclyc sum, is
another tensor which plays a key r{\^o}le in ${\sl G}_2$-structures.
In terms of a Cayley frame, $\ast \varphi$ is locally given by
\begin{gather*}
\ast \varphi = - \sum_{i \in \mathbb{Z}_7} e_{i+2} \wedge e_{i+4}
\wedge e_{i+5} \wedge e_{i+6}.
\end{gather*}
As volume form we fix the form $Vol$ such that $\varphi \wedge
\ast \varphi = 7 Vol = 7 e_0 \wedge \dots \wedge e_6$.

If $\overline{\nabla}^{{\sl G}_2} = \overline{\nabla} + \xi^{{\sl
G}_2}$ is the minimal ${\sl G}_2$-connection, then the Levi-Civita
covariant derivative $\overline{\nabla} \varphi = - \xi^{{\sl
G}_2} \varphi  \in T^* \overline{M} \otimes \lie{g}_2^{\perp}
\subseteq T^* \overline{M} \otimes \Lambda^3 T^* \overline{M}$.
Since $\lie{g}_2^{\perp} = \{ x \lrcorner \ast \varphi \, | \, x
\in T\overline{M} \}$,  the space ${\cal X}$ of covariant
derivatives of $\varphi$ can be described in the way contained in
the following lemma.
\begin{lemma}[\cite{Cabrera:G2struct}] \label{ALDER}
If $ \{ e_0 , \dots , e_6 \} $ is a Cayley frame, then $\alpha \in
{\cal X}$  if and only if
\begin{equation} \label{aijcova}
\alpha  =  \sum_{ i,j \in {\Bbb Z}_{\; 7}} a_{ij} \;  e_{i}
\otimes e_j \lrcorner \ast  \varphi = - \sum_{ i,j \in {\Bbb
Z}_{\; 7}} a_{ij} \;  e_{i} \otimes \ast (e_j \wedge \varphi).
\end{equation}
\end{lemma}

Therefore, if $\overline{\nabla}$ denotes the Levi-Civita
connection of the metric  tensor field $\langle \cdot , \cdot
\rangle$ on $\overline{M}$, the tensor $\overline{\nabla} \varphi$
can be expressed as $\alpha$ in Equation \eqref{aijcova}.
Moreover, the exterior derivative $d \varphi$ and the coderivative
$d^* \varphi$ are given by
  \begin{gather} \nonumber
d \varphi =   - \displaystyle \sum_{i \in {\Bbb Z}_{\;7}} \left(
a_{i+2,i+2} + a_{i+4,i+4} + a_{i+5,i+5} + a_{i+6,i+6} \right) \,
e_{i+2} \wedge e_{i+4} \wedge  e_{i+5} \wedge e_{i+6}
\\ \nonumber
 \; \; \;+ \displaystyle  \sum_{i \in {\Bbb Z}_{\;7}} \left(  a_{i+4,i+5}
+ a_{i+1,i+3} + a_{i+2,i+6}\right) \, e_{i} \wedge e_{i+1} \wedge
e_{i+2} \wedge e_{i+4} \\
\label{diffund}  + \displaystyle  \sum_{i \in {\Bbb Z}_{\;7}}
\left( -a_{i+5,i+4} -a_{i+3,i+1} + a_{i+2,i+6}\right) \, e_{i}
\wedge e_{i+2} \wedge e_{i+3} \wedge e_{i+5} \\
\nonumber  \; \; \; \; \;\, + \displaystyle \sum_{i \in {\Bbb
Z}_{\;7}} \left( a_{i+4,i+5} -a_{i+3,i+1}
-a_{i+6,i+2}\right) \, e_{i} \wedge e_{i+3} \wedge e_{i+4} \wedge e_{i+6} \\
\nonumber  \; \; \; \; \;\;\;\; +  \displaystyle  \sum_{i \in
{\Bbb Z}_{\;7}} \left(- a_{i+5,i+4} + a_{i+1,i+3} -
a_{i+6,i+2}\right)\, e_{i} \wedge e_{i+5} \wedge e_{i+6} \wedge
e_{i+1}, \\
 \nonumber d^*  \varphi =    \sum_{i \in {\Bbb Z}_{\;7}} \left(
a_{i+5,i+4} - a_{i+4,i+5} + a_{i+6,i+2} - a_{i+2,i+6} \right) \,
e_{i+1} \wedge e_{i+3}
\\
\label{coderivada} \; \; \; \; \;\;\;\;   +  \sum_{i \in {\Bbb
Z}_{\;7}} \left( a_{i+3,i+1} - a_{i+1,i+3} + a_{i+6,i+2} -
a_{i+2,i+6}\right) \, e_{i+4}
\wedge e_{i+5}  \\
\nonumber \; \; \; \; \;\;\;\;  \;\;\;\;  \;\;\;\; +  \sum_{i \in
{\Bbb Z}_{\;7}} \left( a_{i+3,i+1} -a_{i+1,i+3} + a_{i+5,i+4} -
a_{i+4,i+5} \right) \,e_{i+2} \wedge e_{i+6}.
\end{gather}

\begin{table}[tp]
  \centering
   {\footnotesize
  \begin{tabular}{ll}
    \toprule
    $ {\cal P}$ & $\overline{r}(\varphi)=0$ \\
    \midrule
${\cal X}_1 = {\cal NP}$ & $\overline{r}(\varphi)=\frac{k}{4} \langle \cdot , \cdot \rangle$, $k$ constant \\
\midrule
 ${\cal X}_2 ={\cal AP} $ & $\overline{r}(\varphi) \in
{\frak
g}_2$ \\
\midrule
${\cal X}_3$ & $\overline{r}(\varphi) \in S_{0}^2 T^{\ast} \overline{M} $  \\
\midrule
${\cal X}_4 = {\cal LCP}$ & $\overline{r}(\varphi) = - \frac{1}{12} p_{\varphi} \lrcorner \varphi$ \\
\midrule ${\cal X}_1 + {\cal X}_2 $ & $\overline{r}(\varphi) +
\overline{r}(\varphi)^t $ is scalar and  $\overline{r}(\varphi) -
\overline{r}(\varphi)^t \in {\frak g}_2$   \\
\midrule ${\cal X}_1 + {\cal X}_3 = {\cal SP} $ &
$\overline{r}(\varphi)$ is symmetric \\
\midrule ${\cal X}_2 + {\cal X}_3 $ & $\overline{r}(\varphi) +
\overline{r}(\varphi)^t \in S_0^2 T^{\ast} \overline{M}$ and
$\overline{r}(\varphi) -
\overline{r}(\varphi)^t \in {\frak g}_2$  \\
\midrule ${\cal X}_1 + {\cal X}_4 = {\cal LCNP} $ &
$\overline{r}(\varphi) + \overline{r}(\varphi)^t $ is scalar and
$\overline{r}(\varphi) -
\overline{r}(\varphi)^t =  - \frac16 p_{\varphi} \lrcorner \varphi$  \\
\midrule
${\cal X}_2 + {\cal X}_4 ={\cal LCAP} $ & $\overline{r}(\varphi)$ is skew-symmetric \\[1mm]
\midrule
 ${\cal X}_3 + {\cal X}_4 $ &  $\overline{r}(\varphi) +
\overline{r}(\varphi)^t \in S_0^2 T^{\ast} \overline{M}$ and
$\overline{r}(\varphi) -
\overline{r}(\varphi)^t = -  \frac16 p_{\varphi} \lrcorner \varphi $      \\
\midrule
${\cal X}_1 + {\cal X}_2 + {\cal X}_3$ & $\overline{r}(\varphi) - \overline{r}(\varphi)^t \in {\frak g}_2 $ \\
\midrule
 ${\cal X}_1 + {\cal X}_2 + {\cal X}_4 $ &
$\overline{r}(\varphi) +
\overline{r}(\varphi)^t $ is scalar \\
\midrule ${\cal X}_1 + {\cal X}_3 + {\cal X}_4 $ &
$\overline{r}(\varphi) - \overline{r}(\varphi)^t =  - \frac16 p_{\varphi} \lrcorner \varphi$  \\
\midrule
 ${\cal X}_2 + {\cal X}_3 + {\cal X}_4 $ &
$\overline{r}(\varphi) +
\overline{r}(\varphi)^t \in S_0^2 T^{\ast} \overline{M}$ \\
\midrule
${\cal X} $ &  no relation \\
    \bottomrule
  \end{tabular}  }
  \caption{Types of $\Lie{G}_2$-structures}
  \label{tab:g2types}
\end{table} \vspace{-2mm}

\noindent Note that one can compute $a_{ij}$ from $d\varphi$ and
$d^* \varphi$ , for all $i,j \in {\Bbb Z}_{\;7}$. As a
consequence, the expression for $\overline{\nabla} \varphi$ can be
rebuilt from  $d\varphi$ and $d^* \varphi$. \vspace{2mm}

 Fern{\'a}ndez and  Gray \cite{Fern-Gray:G2struct} proved that
${\cal X}$, under the action of $\Lie{G}_2$, has four irreducible
components, $ {\cal X} = {\cal X}^{(1)}_1 + {\cal X}^{(14)}_2 +
{\cal X}^{(27)}_3 + {\cal X}^{(7)}_4$,  where the upper index
indicates the corresponding dimension. Thus $\overline{\nabla}
\varphi \in {\cal X}$ has four components giving rise to sixteen
types of ${\sl G}_2$-structure. Conditions in terms of
 $d\varphi$ and $d^* \varphi$ can be given to characterize each type
 of ${\sl G}_2$-structure \cite{Cabrera:G2struct}.
    In \cite{Cabrera:Palermo}, it was
showed that the irreducible $\Lie{G}_2$-summands of ${\cal X}$ can
be alternatively described by the following
$\Lie{G}_2$-equivariant map
\begin{equation}
\begin{array}{rcl}
{\cal X} & \to & T^{\ast} \overline{M} \otimes T^{\ast} \overline{M} \\
\alpha   & \to & \overline{r}(\alpha),
\end{array}
\end{equation}
where $\overline{r}(\alpha)(x,y) = \frac1{4} \langle x \lrcorner
\alpha , y \lrcorner \ast \varphi \rangle$. For all $\alpha \in
{\cal X}$ given by Equation \eqref{aijcova}, we have
\begin{equation}\label{ralfa}
\overline{r}(\alpha) = \displaystyle \sum_{i , j \in {\Bbb
Z}_{\;7}} a_{ij} \, e_{i} \otimes e_{j}.
\end{equation}
Hence  $\overline{r}$ is a $\Lie{G}_2$-isomorphism. For the
covariant two-tensors on $\overline{M}$, we  have the following
decomposition into $\Lie{G}_2$-irreducible components, $ \otimes^2
T^{\ast} \overline{M} = {\Bbb R} + \Lie{S}_0^{2} T^{\ast}
\overline{M} + \lie{g}_2 + \lie{g}^{\perp}_2 $,
 where $\Lie{S}_0^{2} T^{\ast}\overline{M}$ is the space of trace free
symmetric two-tensors, $\lie{g}_2$ is the Lie algebra of
$\Lie{G}_2$ and $\lie{g}_2^{\perp}$ is the orthogonal complement
of $\lie{g}_2$ in $\Lambda^2 T^{\ast} \overline{M}$. Now by the
$\Lie{G}_2$-isomorphism $\overline{r}$, using the Schur's Lemma,
we get ${\cal X}^{(1)}_1 \cong {\Bbb R}$,  ${\cal X}^{(14)}_2
\cong \lie{g}_2$, ${\cal X}^{(27)}_3 \cong \Lie{S}^2_0 T^{\ast}
\overline{M}$ and ${\cal X}^{(7)}_4 \cong \lie{g}_2^{\perp}$. We
recall that $\lie{g}_2$ can be described as consisting of those
skew-symmetric bilinear forms $b$ such that
\begin{equation} \label{algebragedos}
 \sum_{\{j,k \, | \, \varphi(e_i,e_j,e_k)=1 \}}
b (e_j , e_k) = 0,
\end{equation}
for all $i \in {\Bbb Z}_{\;7}$, where $\{ e_0 , \dots , e_6 \}$ is
a Cayley frame.

For sake of simplicity we will write $\overline{r}(\varphi)=
\overline{r}(\overline{\nabla} \varphi)$. Now, by the map
$\overline{r}$, using Schur's Lemma,
 one can  characterize  each type of $\Lie{G}_2$-structure
(\cite{Cabrera:Palermo}). Such characterizations are shown in
 Table~\ref{tab:g2types}.

The vector field $p_{\varphi}$ which occurs in some conditions
contained  in  Table~\ref{tab:g2types} is such that $ pd^* \varphi
= \ast (\ast d \varphi \wedge \varphi) = \langle p_{\varphi},
\cdot \rangle $. In \cite{Cabrera:hipspin}, it is showed that
\begin{equation} \label{pdeltag}
 p d^* \varphi  = - 2 \displaystyle \sum_{i \in {\Bbb Z}_{\;7}}
\sum_{\{j,k \in {\Bbb Z}_{\;7} |  \varphi(e_i,e_j,e_k)=1 \}}
\left( \overline{r}(\varphi) (e_j , e_k) -
\overline{r}(\varphi)(e_k,e_j) \right) \, e_{i}.
\end{equation}

It is well known  that the intrinsic torsion $\xi^{{\sl G}_2}$ of
a ${\sl G}_2$-structure can be computed from $d \varphi$ and $d
\ast \varphi$. Here we will give an explicit description of this
fact.

\begin{theorem} The minimal ${\sl G}_2$-connection is given by
$\widetilde{\nabla}^{\Lie{G}_2} = \overline{\nabla} + \xi^{{\sl
G}_2}$, where $\xi^{{\sl G}_2}$ is defined by
$$
\xi^{{\sl G}_2}_X Y =  - \frac13 \sum_{i \in {\mathbb Z}_7}
\overline{r}(\varphi) (X, e_i) P(e_i,Y),
$$
for all vectors $X,Y$, and $\overline{r}(\varphi)= 1/4 \langle
\overline{\nabla}_{\cdot} \varphi , \cdot \lrcorner \ast \varphi
\rangle$. Moreover, the bilinear form $\overline{r}(\varphi)$ is
expressed in  terms of $d \varphi$ and $\ast \varphi$ by
\begin{equation} \label{rxy}
4 \overline{r}(\varphi)(X,Y) =  \langle X \lrcorner d \varphi, Y
\lrcorner \ast \varphi ,
 \rangle -  \langle Y \lrcorner ( X \wedge
\ast \varphi) , d \varphi \rangle + 2 d^* \varphi (X,Y).
\end{equation}
\end{theorem}
\begin{proof} Equation (\ref{rxy}) can be checked for
$\overline{r}(\varphi)(e_i,e_j)$, using the expressions
(\ref{diffund}) and (\ref{coderivada}) for $d\varphi$ and
$d^{\ast} \varphi$, respectively.

It is also immediate that $\xi^{{\sl G}_2} \in T^* M \otimes
\lie{g}_2^{\perp}$. Finally, it is straightforward to check that
$\widetilde{\nabla}^{\Lie{G}_2} \varphi = 0$. Hence
$\widetilde{\nabla}^{\Lie{G}_2}$ is the minimal
$\Lie{G}_2$-connection.
\end{proof}

\begin{remark}
As a direct consequence of last Theorem, we obtain an expression
for $\overline{\nabla} \varphi$ in terms of $d\varphi$ and
$d^{\ast} \varphi$, i.e., $\overline{\nabla} \varphi = - \xi^{{\sl
G}_2} \varphi$. Also note that the $\Lie{G}_2$-connection
$\widetilde{\nabla}^{\Lie{G}_2}$, particularised   to
$\Lie{G}_2$-structures of type  ${\cal X}_2$, coincides with the
one given by Cleyton and Ivanov in \cite{CleytonIvanov}.
\end{remark}

\section{Orientable hypersurfaces of  manifolds with $\Lie{G}_2$-structure}
\label{sec:orienthyp}
 From
now on, $M$ will be an orientable hypersurface of a
seven-dimensional Riemannian manifold $\overline{M}$ with a
$\Lie{G}_2$-structure and $\iota \, : \, M \to \overline{M}$ will
denote the inclusion map. Associated with the
$\Lie{G}_2$-structure, we have the  metric $\langle \cdot , \cdot
\rangle$, the fundamental three-form $\varphi$ and the two-fold
vector cross product $P$.
\begin{proposition} \label{g2hyp}
Let $\overline{M}$ be a Riemannian manifold with a
$\Lie{G}_2$-structure. If
 $M$ is an orientable hypersurface and $n$ is a unit normal vector
field on $M$, then there is a special almost Hermitian structure
on $M$ defined by the almost complex structure
\begin{equation} \label{phip}
Ix = P(n,x),
\end{equation}
and the complex volume form given by
\begin{gather}
\label{volumencomplejo}
 \Psi  = \cos \theta
\iota^* \varphi   - \sin \theta \iota^* ( n \lrcorner \ast
\varphi) + i \, \left( \sin \theta  \iota^* \varphi + \cos \theta
\iota^* ( n
\lrcorner \ast \varphi) \right), 
\end{gather}
where $\theta$ is a smooth function on $M$.
\end{proposition}

The almost complex structure defined by \eqref{phip} has been
already  considered  by Calabi in \cite{Calabi:ac6manifold} and by
Gray in \cite{Gray:crossproduct}. On each point of $M$ there is a
 Cayley basis with $n$ as first element, i.e., $\{n, e_1, \dots ,
e_6 \}$. The local frame of $M$ given by $\{ e_1, e_2 , e_4 , e_3
, e_6 , e_5 \}$ is an adapted local frame for the special almost
Hermitian structure defined in Proposition \ref{g2hyp} and the
special unitary group $\SU(3)$ can be considered as included in
$\Lie{G}_2$ in the following way
$$
\SU(3) = \{ g \in \Lie{G}_2 \subset \SO(7) \; | \; g.n=n \; \}.
$$

In next lemma we relate the bilinear form $r(\omega)= r(\nabla
\omega) $, defined by Equation \eqref{erresu}, with the bilinear
form $\overline{r}(\varphi)$ and the shape tensor $II$.

\begin{proposition} \label{rrB}
For the $\SU(3)$-structure and  the $\Lie{G}_2$-structure
considered in Proposition \ref{g2hyp}, if $B(x,y) = \langle
II(x,y),n\rangle$ denotes the second fundamental form, $\Lie{Tr}$
denotes trace, $c_{\omega}$ means the contraction by $\omega$ and
$h$ is the length of the mean curvature $H$, i.e., $h = \langle H
, n \rangle$, then
\begin{eqnarray}
  r(\omega)& = &\cos \theta
\left( - I_{(2)} \iota^* \overline{r}(\varphi)  + B \right) - 
\sin \theta \left(   \iota^* \overline{r}(\varphi) + I_{(2)} B
\right), \label{rrB1} \\
  2 I d^* \omega & = &  \iota^*  p d^*
 \varphi  - 2 \overline{r}(\varphi) (n , \iota_{\ast} I \cdot )
 + 2 \overline{r}(\varphi)( \iota_{\ast} I\cdot , n
 ), \label{rrB2}\\
   p d^* \varphi (n) & = & - 2 \cos \theta \,
\Lie{Tr} \left( r(\omega) \right) - 2 \sin \theta \, c_{\omega}
\left( r(\omega) \right)  + 12 h,  \label{rrB3} \\ \qquad Tr
\left( \iota^* \overline{r}(\varphi) \right) & = &  - \sin\theta
\, \Lie{Tr} \left( r(\omega) \right) + \cos\theta \,
c_{\omega}\left(
r(\omega) \right), \label{rrB4} \\
3 I\eta &=  &d \theta  - \overline{r}(\varphi) ( \iota_{\ast}
\cdot , n ). \label{rrB5}
\end{eqnarray}
\end{proposition}
\begin{proof}
On each point of $M$, we consider a  Cayley frame $\{ n, e_{1},
\dots , e_6 \}$ and, using Equation \eqref{ralfa} and Lemma
\ref{ALDER}, obtain
\begin{gather*}
\overline{r}(\varphi) (e_i , e_1) = \overline{a}_{i1} = \left(
\overline{\nabla}_{e_i} \varphi \right) (e_4 , e_6 , n ) = \langle
\left( \overline{\nabla}_{e_i} P \right) (n , e_4) , e_6 \rangle\\
= \langle \left( \nabla_{e_i} I \right) e_4 , e_6 \rangle -
\langle \overline{\nabla}_{e_i} n , P(e_4, e_6)
\rangle\\
= - \left( \nabla_{e_i} \omega \right) (e_4 , e_6) + B(e_i ,
Ie_1).
\end{gather*}
Since we have
$$
2  \left( \nabla_{e_i} \omega \right) (e_4 , e_6) = \langle
\nabla_{e_i} \omega , Ie_1 \lrcorner (\iota^* \varphi) \rangle  =
- \langle \nabla_{e_i} \omega , e_1 \lrcorner \iota^* (n \lrcorner
\ast \varphi) \rangle,
$$
it is not hard to show
\begin{gather*}
\overline{r}(\varphi) (X , IY) =   \frac{1}{2} \langle \nabla_X
\omega , Y \lrcorner (\iota^* \varphi) \rangle - B(X,Y),\\
\overline{r}(\varphi) (X , Y) =   \frac{1}{2} \langle \nabla_X
\omega , Y \lrcorner \iota^*(n \lrcorner \ast \varphi) \rangle +
B(X,IY).
\end{gather*}
for all vectors $X,Y \in T M$. From these two identities, Equation
\eqref{rrB1} follows.

Equation~\eqref{rrB2} is deduced from  Equation \eqref{pdeltag},
taking Equations~\eqref{rrB1} and~\eqref{codeseis} into account.
Equations~\eqref{rrB3} and~\eqref{rrB4}  are derived by computing
$ {\textsl Tr} \left( r(\omega) \right)$ and $c_{\omega} \left(
r(\omega) \right)$, taking Equations~\eqref{pdeltag}
and~\eqref{rrB1}  into account.

For Equation~\eqref{rrB5}, we firstly consider that $\theta$ is
constant and equal to $0$. In such a case, we have $\psi_+ = i^*
\varphi$ and $\psi_- = \iota^* ( n \lrcorner \ast \varphi)$.
Noting  that $d i^* \varphi = i^* d \varphi$ and
 making use of the  explicit expression
for $d \varphi$ given by Equation \eqref{diffund}, we obtain the
following identities
\begin{eqnarray*}
\ast \left( \ast d \, \iota^* \varphi \wedge \iota^* \varphi
\right) & = & I
d^* \omega - 2 \overline{r}(\varphi) (\iota_{\ast} I\cdot , n), \\
\ast \left( \ast d \, \iota^* \varphi \wedge \iota^* (n \lrcorner
\ast \varphi) \right) & = & - d^* \omega - 2 \overline{r}(\varphi)
(\iota_{\ast} \cdot , n).
\end{eqnarray*}
From these identities, taking Theorem~\ref{torsion3w5} into
account, it follows
\begin{eqnarray}  \label{etaprimer}
\ast \left( \ast d \iota^* (n \lrcorner \ast \varphi) \wedge
\iota^*(n \lrcorner \ast \varphi) \right) & = & \ast \left( \ast
\iota^* d
\varphi \wedge \iota^* \varphi \right) \\
\nonumber & = & I d^* \omega - 2
\overline{r}(\varphi) (\iota_{\ast} I\cdot,n),\\
- \ast \left( \ast d \iota^* (n \lrcorner \ast \varphi) \wedge
\iota^* \varphi \right)   & = & \ast \left( \ast d \, \iota^*
\varphi
\wedge \iota^* (n \lrcorner \ast \varphi) \right)  \label{etasegundo} \\
\nonumber &  = &  - d^* \omega - 2 \overline{r}(\varphi) (
\iota_{\ast} \cdot , n).
\end{eqnarray}

On the other hand, making use of Equations \eqref{estrella1} and
\eqref{estrella2}, for all $\mu \in T^* M$, we have
\begin{eqnarray} \label{etatercer}
&&\ast \left( \ast ( \mu \wedge \iota^* \varphi ) \wedge \iota^*
\varphi \right) = \ast \left( \ast ( \mu \wedge \iota^* (n
\lrcorner \ast \varphi) ) \wedge \iota^* (n \lrcorner \ast
\varphi) \right)  = -2 \mu,
\\
 && \ast \left( \ast ( \mu \wedge \iota^* (n \lrcorner \ast \varphi) ) \wedge
 \iota^* \varphi \right) = - \ast \left( \ast ( \mu \wedge \iota^* \varphi )
\wedge \iota^* (n \lrcorner \ast \varphi) \right) =  2 I\mu.
\label{etacuarto}
\end{eqnarray}

For sake of simplicity we have deduced Equations
\eqref{etaprimer},
 \eqref{etasegundo}, \eqref{etatercer} and \eqref{etacuarto}  fixing
  $\Psi = \iota^* \varphi + i \iota^* ( n\lrcorner
\ast \varphi)$. But really, these equations  mainly depend of the
Hodge star operator $\ast$ which is determined by the metric, and
not of what complex volume form $\Psi$ we have previously fixed.
Thus Equations \eqref{etaprimer}, \eqref{etasegundo},
\eqref{etatercer} and \eqref{etacuarto}  are still true when
$\Psi$ is given by \eqref{volumencomplejo}. In such a case, we
would have $\psi_+ = \cos \theta  \, \iota^* \varphi   - \sin
\theta \, \iota^* (n \lrcorner \ast \varphi)$ and need to compute
$\ast \left( \ast d \psi_+ \wedge \psi_+ \right)$ to determine
$\eta$ (see Theorem \ref{torsion3w5}). In fact, taking Equations
\eqref{etaprimer}, \eqref{etasegundo}, \eqref{etatercer} and
\eqref{etacuarto} into account, we would get
$$
6 \eta + I d^* \omega = \ast \left( \ast d \psi_+ \wedge \psi_+
\right)= - 2 I d \theta + I d^* \omega - 2 \overline{r}(\varphi) (
\iota_{\ast} I \cdot , n).
$$
From these identities Equation~\eqref{rrB5} follows.
 \end{proof}

Proposition \ref{rrB},  Table~\ref{tab:g2types} and Proposition
\ref{suseis} give place to results relating the type of
$\Lie{G}_2$-structures on the ambient manifold, type of the
$\SU(3)$-structure on the hypersurface and the shape tensor. In
the following theorem we only mention  the more relevant
consequences  in such a direction.
\begin{theorem} \label{xunogeneral}
Let $\overline{M}$ be a seven-dimensional Riemannian manifold with
a
 $\Lie{G}_2$-structure of type ${\cal X}_1$ such that $d
\varphi = k \ast \varphi$. Let $M$ be an orientable hypersurface
with unitary normal vector field $\,n$. We consider an
$\SU(3)$-structure on $M$  defined as in Proposition \ref{g2hyp}.
Then $M$ is of type  $ {\cal W}^+_1 + {\cal W}_1^- + {\cal W}^+_2
+ {\cal W}_2^- + {\cal W}_3 + {\cal W}_5$ and the conditions
displayed in Table~\ref{tab:xunogeneral} characterise types of
$\SU(3)$-structure on $M$.
\end{theorem}
\begin{proof} Taking Proposition~\ref{rrB} into account, since $\overline{r} ( \varphi )=
k/4 \langle \cdot , \cdot \rangle$  and $pd^*\varphi=0$ (see
Table~\ref{tab:g2types}),
 we have
\begin{gather*}
r(\omega)  = \cos \theta \left( \frac{k}{4}  \omega +  B \right) -
\sin \theta \left(
 \frac{k}{4}  \langle \cdot , \cdot \rangle +  I_{(2)}B \right),
\\
\quad I d^* \omega = 0, \quad 3 I \eta = d \theta.
\end{gather*}
Now, all parts of Theorem are direct consequences of these
equations and Proposition~\ref{suseis}.
\end{proof}
\begin{remark} \begin{enumerate}
\item[(i)] Conditions given in 10th, 11th and 12th lines of
Table~\ref{tab:xunogeneral} has been already shown by Calabi
\cite{Calabi:ac6manifold} for the particular case $\overline{M} =
\mathbb R^7$ considered as the pure imaginary Cayley numbers. Also
for such a case, Gray \cite{Gray:crossproduct} proved that any
orientable hypersurface is of type ${\cal W}_1 +  {\cal W}_2 +
{\cal W}_3 $ as almost Hermitian manifold, and deduced conditions
with correspond with those contained in 2nd and 9th lines of
Table~\ref{tab:xunogeneral}. Here the context is more general and
we give more detailed information with respect to the
$\SU(3)$-structure defined on the hypersurface.
 \item[(ii)] In  Theorem~\ref{xunogeneral}, if we
consider  the $\Lie{G}_2$-structure on $\overline{M}$ of type
$\mathcal P$,  we would obtain the same conclusions but for $k=0$.
So that we will not show a specific table for such a situation.
 \item[(iii)] In Table~\ref{tab:xunogeneral}, it is not given an
 exhaustive list of
 conditions which characterise all possible types. We only has
 written conditions for those more relevant types. However, it is straightforward to
 deduce the remaining characterizations  from conditions contained in the
 table.
\end{enumerate}
\end{remark}

\begin{table}[tbp]
  \centering {\footnotesize
  \begin{tabular}{ll}
    \toprule
    $ {\cal W}^+_1 + {\cal W}_1^- + {\cal W}^+_2 + {\cal W}_2^- + {\cal W}_3 $ &  $\theta$ is  constant  \\
    \midrule
    ${\cal W}^+_1 + {\cal W}_1^- + {\cal W}^+_2 + {\cal W}_2^-
    + {\cal W}_5$  &  $IB = B$\\
     \midrule
     ${\cal W}^+_1 + {\cal W}_1^- + {\cal W}^+_2 +  {\cal W}_3 + {\cal W}_5$  &
     $\sin \theta \, (1+I)B = 2 h \sin \theta \langle \cdot , \cdot
     \rangle$ \\
      \midrule
      ${\cal W}^+_1 + {\cal W}_1^-
      + {\cal W}^-_2 +  {\cal W}_3 + {\cal W}_5$ &
      $\cos \theta \, (1+I)B = 2 h \cos \theta \langle \cdot ,
      \cdot\rangle$\\
      \midrule
      ${\cal W}^+_1 + {\cal W}_2^+ +
       {\cal W}^-_2 +  {\cal W}_3 + {\cal W}_5$ &   $- 4 h
       \sin \theta= k \cos \theta$ \\
      \midrule
      ${\cal W}^-_1 + {\cal W}_2^+ + {\cal W}^-_2 +  {\cal W}_3 + {\cal W}_5$
      & $4 h \cos \theta =  k \sin \theta$ \\
      \midrule
      ${\cal W}_1^+ + {\cal W}^-_1 +  {\cal W}_3 + {\cal W}_5$ & $(1+I)B =
      2h
      \langle \cdot , \cdot \rangle$\\
      \midrule
       ${\cal W}_2^+ + {\cal W}^-_2 + {\cal W}_3 + {\cal W}_5$ & $M$ is a minimal variety and $\overline{M}$ is of
        type ${\cal P}$\\
       \midrule
       ${\cal W}^+_1+ {\cal W}^-_1 + {\cal W}_5$ & $M$ is totally umbilic\\
       \midrule
        ${\cal W}^+_2 + {\cal W}_2^- + {\cal W}_5$ & $IB=B$, $M$ is a minimal
         variety and  $\overline{M}$ is of type ${\cal P}$\\
       \midrule
       ${\cal W}_3 + {\cal W}_5$ & $I B=-B$ and
       $\overline{M}$ is of type ${\cal P}$\\
       \midrule
       ${\cal W}_5$ & $M$ is totally geodesic and  $\overline{M}$ is of type ${\cal
       P}$\\
       \midrule
       $\{ 0 \}$ & $M$ is totally geodesic, $\theta$ is constant and \\
       & $\overline{M}$ is of
       type ${\cal P}$\\
    \bottomrule
  \end{tabular}  }
  \caption{$\overline{M}$ of type ${\cal X}_1$}
  \label{tab:xunogeneral}
\end{table}

For most of the following theorems, their proves are deduced using
analog arguments as in the proof of  Theorem~\ref{xunogeneral}.
For such a reason, we give some of  those  theorems without an
explicit proof.

Two special particular cases of Theorem~\ref{xunogeneral} for
$\theta$ constant are respectively  when $\theta = 0$  and $
\theta = \pi /2$.
\begin{theorem} \label{xunopimedio}
Let $\overline{M}$ be a seven-dimensional Riemannian manifold with
a
 $\Lie{G}_2$-structure of type ${\cal X}_1$ such that $d
\varphi = k \ast \varphi$. Let $M$ be an orientable hypersurface
with unitary normal vector field $\,n$. We consider an
$\SU(3)$-structure on $M$  defined as in Proposition~\ref{g2hyp},
taking  $\theta= 0(\theta =\pi /2)$ constant. Then $M$ is of type
$ {\cal W}^{+(-)}_1 + {\cal W}_1^{-(+)} + {\cal W}^{+(-)}_2 +
{\cal W}_3 $ and the conditions displayed in
Table~\ref{tab:xunopimedio} characterise types of
$\SU(3)$-structure on $M$.
\end{theorem}

A particular case of Theorem~\ref{xunopimedio} is when the
$\Lie{G}_2$-structure on $\overline{M}$ is of type ${\cal P}$. In
such a situation we have the following result.

\begin{theorem} \label{xparalelopimedio}
Let $\overline{M}$ be a seven-dimensional Riemannian manifold with
a
 $\Lie{G}_2$-structure of type ${\cal P}$. Let $M$ be an orientable hypersurface
with unitary normal vector field $\,n$. We consider an
$\SU(3)$-structure on $M$  defined as in Proposition~\ref{g2hyp},
taking  $\theta= 0(\theta= \pi / 2)$ constant. Then $M$ is of type
$ {\cal W}^{+(-)}_1 +  {\cal W}^{+(-)}_2 +  {\cal W}_3 $ and the
conditions displayed in Table~\ref{tab:xparalelopimedio}
characterise  types of $\SU(3)$-structure on $M$.
\end{theorem}

\begin{table}[tbp]
  \centering {\footnotesize
  \begin{tabular}{ll}
    \toprule
    ${\cal W}^{+(-)}_1 + {\cal W}_1^{-(+)} + {\cal W}^{+(-)}_2$     &  $IB = B$\\
      \midrule
      ${\cal W}^{+(-)}_1 + {\cal W}_1^{-(+)}
       +  {\cal W}_3$  &
      $ (1+I)B = 2 h \langle \cdot , \cdot\rangle$\\
      \midrule
      ${\cal W}^{+(-)}_1 + {\cal W}_2^{+(-)} +
        {\cal W}_3 $ &  $\overline{M}$ is of
       type ${\cal P}$ \\
      \midrule
      ${\cal W}^{-(+)}_1 + {\cal W}_2^{+(-)}  + {\cal W}_3 $
      & $M$ is a minimal variety  \\
      \midrule
       ${\cal W}^{+(-)}_1+ {\cal W}^{-(+)}_1 $ & $M$ is totally umbilic\\
       \midrule
       ${\cal W}^{+(-)}_1+ {\cal W}^{+(-)}_2 $ &$IB = B$ and  $\overline{M}$ is of
       type ${\cal P}$ \\
        \midrule
       ${\cal W}^{-(+)}_1+ {\cal W}^{+(-)}_2 $ & $IB=B$ and $M$ is a minimal variety  \\
        \midrule
       ${\cal W}^{+(-)}_1+ {\cal W}_3 $ &  $ (1+I)B = 2 h \langle \cdot , \cdot\rangle$ and $\overline{M}$ is of
       type ${\cal P}$  \\
       \midrule
       ${\cal W}^{-(+)}_1+ {\cal W}_3 $ &  $IB=-B$  \\
      \midrule
       ${\cal W}_2^{+(-)} + {\cal W}_3$ & $M$ is a minimal variety and $\overline{M}$ is of
        type ${\cal P}$\\
       \midrule
       ${\cal W}_3 $ & $I B=-B$ and
       $\overline{M}$ is of type ${\cal P}$\\
       \midrule
        ${\cal W}^{+(-)}_2 $ & $IB=B$, $M$ is a minimal
         variety and  $\overline{M}$ is of type ${\cal P}$\\
         \midrule
       ${\cal W}^{-(+)}_1 $ & $M$ is totally geodesic\\
        \midrule
        ${\cal W}^{+(-)}_1 $ & $M$ is totally umbilic
         variety and  $\overline{M}$ is of type ${\cal P}$\\
       \midrule
       $\{ 0 \}$ & $M$ is totally geodesic and $\overline{M}$ is of
       type ${\cal P}$\\
    \bottomrule
  \end{tabular}  }
  \caption{$\overline{M}$ of type ${\cal X}_1$ and $\theta=0(\pi/2)$ constant}
  \label{tab:xunopimedio}
\end{table}

\begin{table}[tbp]
  \centering {\footnotesize
  \begin{tabular}{ll}
    \toprule
    ${\cal W}^{+(-)}_1 + {\cal W}^{+(-)}_2$     &  $IB = B$\\
      \midrule
      ${\cal W}^{+(-)}_1 +   {\cal W}_3$  &
      $ (1+I)B = 2 h \langle \cdot , \cdot\rangle$\\
      \midrule
      $ {\cal W}_2^{+(-)}  + {\cal W}_3 $
      & $M$ is a minimal variety  \\
       \midrule
       ${\cal W}_3 $ & $I B=-B$ \\
       \midrule
        ${\cal W}^{+(-)}_2 $ & $IB=B$ and  $M$ is a minimal
         variety \\
        \midrule
        ${\cal W}^{+(-)}_1 $ & $M$ is totally umbilic \\
       \midrule
       $\{ 0 \}$ & $M$ is totally geodesic \\
    \bottomrule
  \end{tabular}  }
  \caption{$\overline{M}$ of type ${\cal P}$ and $\theta=0(\theta= \pi / 2 ) $ constant}
  \label{tab:xparalelopimedio}
\end{table}
Next, we describe the situation when the $\Lie{G}_2$-structure is
almost parallel.
\begin{theorem} \label{xdosgeneral}
Let $\overline{M}$ be a seven-dimensional Riemannian manifold with
a
 $\Lie{G}_2$-structure of type ${\cal X}_2$  . Let $M$ be an
orientable hypersurface with unitary normal vector field $n$.  We
consider an $\SU(3)$-structure on $M$  defined as in
Proposition~\ref{g2hyp}.Then we have $ d^* \omega = 6 I \eta - 2 d
\theta  = - 2 \overline{r} (\varphi)( \iota_* \cdot , n )$  and
the conditions displayed  in Table~\ref{tab:xdosgeneral}
characterise
 types of  $\SU(3)$-structure on $M$.
\end{theorem}
\begin{table}[tbp]
  \centering {\footnotesize
  \begin{tabular}{ll}
    \toprule
${\cal W}^+_1 +  {\cal W}^-_1 + {\cal W}^{+}_2 + {\cal W}^{-}_2 +
{\cal W}_3 + {\cal W}_4 $
         &  $d \theta =  2 \overline{r} (\varphi)( \iota_* \cdot , n) $ \\
             \midrule
    ${\cal W}^{+}_1 + {\cal W}^{-}_1 + {\cal W}^{+}_2 +  {\cal W}^{-}_2 + {\cal W}_3 + {\cal W}_5 $    &
    $\overline{r}(\varphi) (  \iota_* \cdot , n )=0$ \\
  \midrule
${\cal W}^+_1 +  {\cal W}^-_1 + {\cal W}^{+}_2 + {\cal W}^{-}_2 +
{\cal W}_4 + {\cal W}_5 $
         &  $IB = B$ \\
          \midrule
    ${\cal W}^{+}_1 + {\cal W}^{-}_1 + {\cal W}^{+}_2 +  {\cal W}_3 + {\cal W}_4 + {\cal W}_5 $    &
    $\begin{array}{l} \sin \theta \left( (I_{(1)} - I_{(2)}) \iota^* \overline{r}(\varphi) +
    (1+I)B\right) \\
    \hspace{4.5cm} = 2 h \sin \theta \langle \cdot , \cdot \rangle
    \end{array}$ \\
    \midrule
    ${\cal W}^{+}_1 + {\cal W}^{-}_1 + {\cal W}^{-}_2 +  {\cal W}_3 + {\cal W}_4 + {\cal W}_5 $    &
    $\begin{array}{l}
    \cos \theta \left( (I_{(1)} - I_{(2)}) \iota^* \overline{r}(\varphi) +
    (1+I)B\right)\\
    \hspace{4.5cm}  = 2 h \cos \theta \langle \cdot , \cdot \rangle \end{array}$ \\
          \midrule
    ${\cal W}^+_1 + {\cal W}^{+}_2 + {\cal W}^{-}_2 + {\cal W}_3 + {\cal W}_4 + {\cal W}_5$
         &  $h \sin \theta =0$ \\
    \midrule
    ${\cal W}^-_1 + {\cal W}^{+}_2 + {\cal W}^{-}_2 + {\cal W}_3 + {\cal W}_4 + {\cal W}_5$
         &  $h \cos \theta =0$ \\
    \midrule

    ${\cal W}^{+}_1 + {\cal W}^{-}_1 + {\cal W}^{+}_2 +  {\cal W}^{-}_2 + {\cal W}_3  $    &
    $\overline{r}(\varphi) (  \iota_* \cdot , n )=0$ and $\theta$ is constant \\
    \midrule
    ${\cal W}^{+}_1 + {\cal W}^{-}_1 +  {\cal W}_3 + {\cal W}_4 + {\cal W}_5 $    &
    $  (I_{(1)} - I_{(2)}) \iota^* \overline{r}(\varphi) + (1+I)B
    = 2 h \langle \cdot , \cdot \rangle$ \\
    \midrule
    ${\cal W}^{+}_2 + {\cal W}^{-}_2 + {\cal W}_3 + {\cal W}_4 + {\cal W}_5$     &  $M$ is a minimal variety\\
    \bottomrule
  \end{tabular}  }
  \caption{$\overline{M}$ of type ${\cal X}_2$}
  \label{tab:xdosgeneral}
\end{table}
\begin{proof}
Taking conditions of Table~\ref{tab:g2types}  into account, it
follows that, for a $\Lie{G}_2$-structure of type ${\cal X}_2$, we
have
\begin{gather*}
c_{\omega} \left( \iota^* \overline{r}(\varphi) \right) =  0,
\qquad  p d^{\ast}  \varphi  =  0, \qquad
\overline{r}(\varphi)(X,Y) = -
\overline{r}(\varphi)(Y,X),  \\
\overline{r}(\varphi)(X,Y)- \overline{r}(\varphi) (IX,IY) =
\overline{r}(\varphi) ( n , P(n, P(X,Y))),
\end{gather*}
for all vector fields $X,Y$ tangent  to $M$. Now,  Theorem follows
from these equations, Proposition~\ref{rrB} and
Proposition~\ref{suseis}.
\end{proof}
Note that Table~\ref{tab:xdosgeneral} only contains conditions for
some types of $\SU(3)$-structure on $M$. However, conditions for
remaining types can be easily  derived from those which are given
in the mentioned table. Next we give the result corresponding to
Theorem~\ref{xdosgeneral} when $\theta=0(\theta=\pi/2)$ constant.
\begin{theorem} \label{xdospimedio}
Let $\overline{M}$ be a seven-dimensional Riemannian manifold with
a
 $\Lie{G}_2$-structure of type ${\cal X}_2$  . Let $M$ be an
orientable hypersurface with unitary normal vector field $n$.  We
consider an $\SU(3)$-structure on $M$  defined as  in
Proposition~\ref{g2hyp}, taking  $\theta= 0(\theta= \pi / 2)$
constant. Then $M$ is of type  ${\cal W}^{+(-)}_1 + {\cal
W}^{+(-)}_2 + {\cal W}_3 + {\cal W}_4 + {\cal W}_5 $, $\; d^*
\omega = 6 I \eta = - 2 \overline{r} (\varphi)( \iota_* \cdot , n
)$ and the conditions displayed  in Table~\ref{tab:xdospimedio}
characterise types of $\SU(3)$-structure on $M$.
\end{theorem}
\begin{table}[tbp]
  \centering {\footnotesize
  \begin{tabular}{ll}
    \toprule
${\cal W}^{+(-)}_1 + {\cal W}^{+(-)}_2 +  {\cal W}_3 $
         &  $ \overline{r} (\varphi)( \iota_* \cdot , n ) =0 $ \\
             \midrule
${\cal W}^{+(-)}_1  + {\cal W}^{+(-)}_2 + {\cal W}_4 + {\cal W}_5$
         &  $IB = B$ \\
    \midrule
    ${\cal W}^{+(-)}_1  +   {\cal W}_3 + {\cal W}_4 + {\cal W}_5 $    &
    $(I_{(1)} - I_{(2)}) \iota^* \overline{r}(\varphi) + (1+I)B    = 2 h \langle \cdot , \cdot \rangle$ \\
          \midrule
    $ {\cal W}^{+(-)}_2 +  {\cal W}_3 + {\cal W}_4 + {\cal W}_5$
         &  $h$ is a minimal variety \\
    \bottomrule
  \end{tabular}  }
  \caption{$\overline{M}$ of type ${\cal X}_2$ and $\theta=0(\theta= \pi / 2) $ constant}
  \label{tab:xdospimedio}
\end{table}
Now, we will see the  corresponding theorem to
$\Lie{G}_2$-structures of type ${\cal X}_3$.
\begin{theorem} \label{xtresgeneral}
Let $\overline{M}$ be a seven-dimensional Riemannian manifold with
a
 $\Lie{G}_2$-structure of type ${\cal X}_3$  . Let $M$ be an
orientable hypersurface with unitary normal vector field $n$.
 We consider an $\SU(3)$-structure on $M$  defined as  in Proposition~\ref{g2hyp}.
Then $M$ is of type ${\cal W}^{+}_1 + {\cal W}^{-}_1 + {\cal
W}^{+}_2 +  {\cal W}^{-}_2 + {\cal W}_3 + {\cal W}_5 $ and the
conditions displayed in Table~\ref{tab:xtresgeneral} characterise
others types of $\SU(3)$-structure on $M$.
\end{theorem}

As a consequence of last Theorem, we analyse  the situation when
$\theta=0 (\theta= \pi /2 )$ constant.
\begin{theorem} \label{xtrespimedio}
Let $\overline{M}$ be a seven-dimensional Riemannian manifold with
a
 $\Lie{G}_2$-structure of type ${\cal X}_3$. Let $M$ be an
orientable hypersurface with unitary normal vector field $n$.
 We consider an $\SU(3)$-structure on $M$  defined as  in Proposition~\ref{g2hyp},
taking  $\theta= 0(\theta= \pi / 2 )$ constant.  Then $M$ is of
type ${\cal W}^{+}_1 + {\cal W}^{-}_1 + {\cal W}^{+}_2 + {\cal
W}^{-}_2 + {\cal W}_3 + {\cal W}_5 $ and the conditions displayed
in Table~\ref{tab:xtrespimedio} characterise types of
$\SU(3)$-structure on $M$.
\end{theorem}

Finally, we pay attention to locally conformal parallel
$\Lie{G}_2$-structures.
\begin{theorem} \label{xcuatrogeneral}
Let $\overline{M}$ be a seven-dimensional Riemannian manifold with
a
 $\Lie{G}_2$-structure of type ${\cal X}_4$  . Let $M$ be an
orientable hypersurface with unitary normal vector field $n$.
  We consider an $\SU(3)$-structure on $M$  defined as in Proposition~\ref{g2hyp}.
Then $3 I d^* \omega =  \iota^* p d^* \varphi =   12 Id \theta +
36  \eta $  and the diverse types of $\SU(3)$-structure on $M$ are
characterised by the conditions displayed in
Table~\ref{tab:xcuatrogeneral}.
\end{theorem}
\begin{proof}
Taking conditions of Table~\ref{tab:g2types}  into account, it
follows that, for a $\Lie{G}_2$-structure of type ${\cal X}_4$, we
have
\begin{gather*}
 2 c_{\omega} \left( r ( \omega ) \right) = \sin \theta (12 h -pd^* \varphi(n)),  \qquad  2
\Lie{Tr}  \left( r(\omega) \right) =  \cos \theta (12 h -pd^* \varphi(n)), \\
 pd^* \varphi(\iota_* I \cdot)= 12  \overline{r} (\varphi)(
\iota_* \cdot, n) =  \varphi( p_{\varphi} , n , \cdot ),
\\
\overline{r}(\varphi)(X,Y) = - \frac{1}{12} pd^*\varphi
(P(X,Y)),\\
pd^* \varphi (P(X,IY)) - pd^* \varphi (P(IX,Y) = 2 pd^* \varphi
(n) \langle X , Y \rangle,
\end{gather*}
for all vector fields $X,Y$ tangent  to $M$. Now,  Theorem follows
from these equations, Proposition~\ref{rrB} and
Proposition~\ref{suseis}.
\end{proof}
\begin{theorem} \label{xcuatropimedio}
Let $\overline{M}$ be a seven-dimensional Riemannian manifold with
a
 $\Lie{G}_2$-structure of type ${\cal X}_4$  . Let $M$ be an
orientable hypersurface with unitary normal vector field $n$.
  We consider an $\SU(3)$-structure on $M$  defined as in Proposition~\ref{g2hyp},
taking  $\theta= 0(\theta= \pi /2 )$ constant. Then $3 I d^*
\omega =  \iota^* p d^* \varphi =   36  \eta $ , $M$ is of type
${\cal W}^{+(-)}_1 + {\cal W}^{+(-)}_2 + {\cal W}_3 + {\cal W}_4
 + {\cal W}_5 $ and the diverse types
of $\SU(3)$-structure on $M$ are characterised by the conditions
displayed in Table~\ref{tab:xcuatropimedio}.
\end{theorem}

Finally, it is of some interest to consider  the particular
situation  such that  the vector field $p_{\varphi}$ is tangent to
$M$.
\begin{corollary} \label{xcuatrocor}
In the same conditions as in Theorem~\ref{xcuatrogeneral}, but
with $p_{\varphi} $ tangent to $M$. Then the diverse types of
$\SU(3)$-structure on $M$ are characterized   by the conditions
given in Table~\ref{tab:xcuatrocor}.
\end{corollary}
\begin{table}[tbp]
  \centering {\footnotesize
  \begin{tabular}{ll}
    \toprule
${\cal W}^+_1 +  {\cal W}^-_1 + {\cal W}^{+}_2 + {\cal W}^{-}_2 +
{\cal W}_3  $
         &  $d \theta =   \overline{r} (\varphi)( \iota_* \cdot , n ) $ \\
       \midrule
    ${\cal W}^{+}_1 + {\cal W}^{-}_1 + {\cal W}^{+}_2 +  {\cal W}^{-}_2 + {\cal W}_5 $    &
    $(I_{(1)} + I_{(2)} ) \overline{r}(\varphi)= (1-I)B$ \\
     \midrule
${\cal W}^+_1 +  {\cal W}^-_1 + {\cal W}^{+}_2 + {\cal W}_3 +
{\cal W}_5 $
         &  $  3 \cos \theta (1+I) \overline{r}(\varphi) + 3 \sin \theta
         (1+I)B$ \\
         & \hspace{2cm} $ = ( 6 h \sin \theta
           - \cos \theta  \overline{r}(\varphi)(n,n)) \langle \cdot , \cdot \rangle$ \\
  \midrule
${\cal W}^+_1 +  {\cal W}^-_1 + {\cal W}^{-}_2 + {\cal W}_3
 + {\cal W}_5 $
         &  $ - 3 \sin \theta (1+I) \overline{r}(\varphi) + 3 \cos \theta
         (1+I)B$ \\
         &\hspace{2cm} $ =  (6 h  \cos \theta  + \sin \theta  \overline{r}(\varphi)(n,n))
          \langle \cdot , \cdot \rangle$ \\
          \midrule
          ${\cal W}^+_1 +  {\cal W}^+_2 + {\cal W}^{-}_2 + {\cal W}_3 + {\cal W}_5 $
         &  $ \cos \theta  \overline{r}(\varphi) (n,n)  =  6 h \sin \theta$ \\
          \midrule
            ${\cal W}^-_1 +  {\cal W}^+_2 + {\cal W}^{-}_2 + {\cal W}_3 + {\cal W}_5 $
         &  $ \sin \theta  \overline{r}(\varphi) (n,n)  =  - 6 h \cos \theta$ \\
          \midrule
          ${\cal W}^+_1 +  {\cal W}^-_1 + {\cal W}_3 + {\cal W}_5 $
         &  $ 3 (1+I) \iota^* \overline{r}(\varphi) =
         -  \overline{r}(\varphi)(n,n) \langle \cdot , \cdot \rangle$ and \\
         &           $(1+I)B = 2 h  \langle \cdot , \cdot \rangle$ \\
         \midrule
    ${\cal W}^{+}_2 + {\cal W}^{-}_2 + {\cal W}_3 +  {\cal W}_5  $    &
    $\overline{r}(\varphi) (n, n)=0$ and $M$ is a minimal  variety \\
    \bottomrule
  \end{tabular}  }
  \caption{$\overline{M}$ of type ${\cal X}_3$}
  \label{tab:xtresgeneral}
\end{table}

\begin{table}[tbp]
  \centering {\footnotesize
  \begin{tabular}{ll}
    \toprule
${\cal W}^+_1 +  {\cal W}^-_1 + {\cal W}^{+}_2 + {\cal W}^{-}_2 +
{\cal W}_3  $
         &  $\overline{r} (\varphi)( \iota_* \cdot , n ) =0 $\\
         & (idem) \\
       \midrule
    ${\cal W}^{+}_1 + {\cal W}^{-}_1 + {\cal W}^{+}_2 +  {\cal W}^{-}_2 + {\cal W}_5 $    &
    $(I_{(1)} + I_{(2)} ) \overline{r}(\varphi)= (1-I)B$ \\
    & (idem) \\
     \midrule
${\cal W}^+_1 +  {\cal W}^-_1 + {\cal W}^{+}_2 + {\cal W}_3 +
{\cal W}_5 $
         &  $ 3  (1+I) \overline{r}(\varphi) =  -  \overline{r}(\varphi)(n,n)) \langle \cdot , \cdot \rangle$ \\
         & ($          (1+I)B =  2 h  \langle \cdot , \cdot
         \rangle$) \\
  \midrule
${\cal W}^+_1 +  {\cal W}^-_1 + {\cal W}^{-}_2 + {\cal W}_3 +
{\cal W}_5 $
         &  $          (1+I)B =  2 h  \langle \cdot , \cdot \rangle$ \\
          & ($ 3  (1+I) \overline{r}(\varphi) =  -  \overline{r}(\varphi)(n,n)) \langle \cdot , \cdot
         \rangle$) \\
          \midrule
          ${\cal W}^+_1 +  {\cal W}^+_2 + {\cal W}^{-}_2 + {\cal W}_3 + {\cal W}_5 $
         &  $   \overline{r}(\varphi) (n,n)  = 0$ \\
              &  ($M$ is a minimal variety) \\
          \midrule
            ${\cal W}^-_1 +  {\cal W}^+_2 + {\cal W}^{-}_2 + {\cal W}_3 + {\cal W}_5 $
         &  $M$ is a minimal variety \\
           &  ($   \overline{r}(\varphi) (n,n)  = 0$) \\
    \bottomrule
  \end{tabular}  }
  \caption{$\overline{M}$ of type ${\cal X}_3$ and $\theta= 0(\theta = \pi / 2)$ constant}
  \label{tab:xtrespimedio}
\end{table}

\begin{table}[tbp]
  \centering {\footnotesize
  \begin{tabular}{ll}
    \toprule
${\cal W}^+_1 +  {\cal W}^-_1 + {\cal W}^{+}_2 + {\cal W}_2^- +
{\cal W}_3 + {\cal W}_4$
         &  $\iota^* p d^* \varphi =   12 Id \theta  $ \\
       \midrule
    ${\cal W}^{+}_1 + {\cal W}^{-}_1 + {\cal W}^{+}_2 + {\cal W}_2^- +  {\cal W}_3 + {\cal W}_5 $    &
    $p_{\varphi}$ is normal to $M$  \\
     \midrule
${\cal W}^+_1 +  {\cal W}^-_1 + {\cal W}^{+}_2 + {\cal W}_2^- +
{\cal W}_4 + {\cal W}_5 $
         &  $ IB=B$ \\
          \midrule
${\cal W}^+_1 +  {\cal W}^-_1 + {\cal W}_2^+  + {\cal W}_3 + {\cal
W}_4
 + {\cal W}_5 $
         & $   \sin \theta \; 2h \langle \cdot
         , \cdot \rangle = \sin \theta (1+I) B
           $ \\
  \midrule
${\cal W}^+_1 +  {\cal W}^-_1 + {\cal W}_2^-  + {\cal W}_3 + {\cal
W}_4
 + {\cal W}_5 $
         & $   \cos \theta \; 2
         h \langle \cdot
         , \cdot \rangle = \cos \theta (1+I) B
           $ \\
          \midrule
          ${\cal W}^+_1 +  {\cal W}^+_2 + {\cal W}_2^-  + {\cal W}_3 + {\cal
W}_4
 + {\cal W}_5 $
         &  $   \sin \theta  (p d^* \varphi(n)  - 12 h)  = 0$ \\
          \midrule
            ${\cal W}^-_1 +  {\cal W}^+_2 + {\cal W}_2^-  + {\cal W}_3 + {\cal
W}_4
 + {\cal W}_5 $
         &  $   \cos \theta  (p d^* \varphi(n) - 12 h)  = 0$ \\
          \midrule
            ${\cal W}^+_1 +  {\cal W}^-_1 +  {\cal W}_3 + {\cal W}_4 + {\cal W}_5 $
         & $(1+I) B = 2 h \langle \cdot , \cdot \rangle
           $ \\
          \midrule
          $  {\cal W}^+_2 + {\cal W}^{-}_2 + {\cal W}_3 + {\cal W}_4 + {\cal W}_5 $
         &  $ p d^* \varphi(n) = 12 h$ \\
          \midrule
          ${\cal W}^+_1 +  {\cal W}^-_1 + {\cal W}_4 + {\cal W}_5 $
         &  $M$ is totally umbilic \\
    \bottomrule
  \end{tabular}  }
  \caption{$\overline{M}$ of type ${\cal X}_4$}
  \label{tab:xcuatrogeneral}
\end{table}
\begin{table}[tbp]
  \centering {\footnotesize
  \begin{tabular}{ll}
    \toprule
    ${\cal W}^{+(-)}_1 +  {\cal W}^{+(-)}_2 + {\cal W}_3 $    &
    $p_{\varphi}$ is normal to $M$  \\
     \midrule
${\cal W}^{+(-)}_1 +  {\cal W}^{+(-)}_2 + {\cal W}_4 + {\cal W}_5
$
         &  $ IB=B$ \\
  \midrule
${\cal W}^{+(-)}_1 +   {\cal W}_3 + {\cal W}_4
 + {\cal W}_5 $
         & $ (1+I) B = 2 h  \langle \cdot
         , \cdot \rangle $ \\
          \midrule
          $  {\cal W}^{+(-)}_2 +  {\cal W}_3 + {\cal W}_4 + {\cal W}_5 $
         &  $ p d^* \varphi(n) =  12 h$ \\
          \midrule
          ${\cal W}^{+(-)}_1 +  {\cal W}_4 + {\cal W}_5 $
         &  $M$ is totally umbilic \\
    \bottomrule
  \end{tabular}  }
  \caption{$\overline{M}$ of type ${\cal X}_4$ and $\theta= 0 (\theta = \pi /2)$ constant}
  \label{tab:xcuatropimedio}
\end{table}

\begin{table}[tbp]
  \centering {\footnotesize
  \begin{tabular}{ll}
    \toprule
${\cal W}^+_1 +  {\cal W}^-_1 + {\cal W}^{+}_2 + {\cal W}_2^- +
{\cal W}_3 + {\cal W}_4$
         &  $\iota^* p d^* \varphi =   12 Id \theta $ \\
     \midrule
${\cal W}^+_1 +  {\cal W}^-_1 + {\cal W}^{+}_2 + {\cal W}_2^- +
{\cal W}_4 + {\cal W}_5 $
         &  $ IB=B$ \\
          \midrule
${\cal W}^+_1 +  {\cal W}^-_1 + {\cal W}_2^+  + {\cal W}_3 + {\cal
W}_4
 + {\cal W}_5 $
         & $  \sin \theta  \; 2   h \langle \cdot
         , \cdot \rangle   = \sin \theta \;(1+I) B
           $ \\
  \midrule
${\cal W}^+_1 +  {\cal W}^-_1 + {\cal W}_2^-  + {\cal W}_3 + {\cal
W}_4
 + {\cal W}_5 $
         &  $  \cos \theta \; 2
         h \langle \cdot
         , \cdot \rangle  = \cos \theta \;(1+I) B )
           $ \\
          \midrule
          ${\cal W}^+_1 +  {\cal W}^+_2 + {\cal W}_2^-  + {\cal W}_3 + {\cal
W}_4
 + {\cal W}_5 $
         &  $  h \sin \theta    = 0$ \\
          \midrule
            ${\cal W}^-_1 +  {\cal W}^+_2 + {\cal W}_2^-  + {\cal W}_3 + {\cal
W}_4
 + {\cal W}_5 $
         &  $  h  \cos \theta   = 0$ \\
          \midrule
            ${\cal W}^+_1 +  {\cal W}^-_1 +  {\cal W}_3 + {\cal W}_4 + {\cal W}_5 $
         &  $ (1+I) B =  2         h \langle \cdot
         , \cdot \rangle
           $ \\
          \midrule
          $  {\cal W}^+_2 + {\cal W}^{-}_2 + {\cal W}_3 + {\cal W}_4 + {\cal W}_5 $
         &  $M$ is a minimal variety \\
          \midrule
          ${\cal W}^+_1 +  {\cal W}^-_1 + {\cal W}_4 + {\cal W}_5 $
         &  $M$ is totally umbilic \\
         \midrule
          ${\cal W}_4 + {\cal W}_5 $
         &  $M$ is totally geodesic \\
    \bottomrule
  \end{tabular}  }
  \caption{$\overline{M}$ of type ${\cal X}_4$ and $p_{\varphi}$ tangent to $M$ }
  \label{tab:xcuatrocor}
\end{table}

Fern{\'a}ndez \& Gray showed that there are at most sixteen types of
seven-dimensional manifolds with $\Lie{G}_2$-structure
\cite{Fern-Gray:G2struct}. 'At most', because topological
conditions can do impossible the existence of
 some particular type. In fact, the type ${\cal W}_1 + {\cal W}_2 -
 ({\cal W}_1 \cup {\cal W}_2)$ can not be found on a
  connected manifold \cite{Cabrera:G2struct}.
 Examples for the remaining types were
 shown in
 \cite{Fern-Gray:G2struct,Fernandez:G2cal,Cabrera:G2struct,CMS}.
 Thus, there are really fifteen types of $\Lie{G}_2$-structure.
 For each type  of $\Lie{G}_2$-structure, theorems and results of the sort here exposed can be
 given. Their proves would make reiterated use of
Proposition~\ref{rrB}, Table~\ref{tab:g2types} and
Proposition~\ref{covaseis}. For sake of brevity, we only  give the
results presented until this point and the following one that we
consider of some special interest.
\begin{theorem} \label{xunotresgeneral}
Let us consider a seven-dimensional Riemannian manifold
$\overline{M}$ equipped with a
 $\Lie{G}_2$-structure,  an
orientable hypersurface $M$ of $\overline{M}$ with unitary normal
vector field $n$, and an $\SU(3)$-structure on $M$ defined as in
Proposition~\ref{g2hyp}.
\begin{enumerate}
    \item[{\rm (i)}] If the $\Lie{G}_2$-structure on $\overline{M}$ is
    of type ${\cal X}_1 + {\cal X}_3$, then the $\SU(3)$-structure on
    $M$ is of type $ {\cal W}^+_1 + {\cal W}_1^- + {\cal W}^+_2 +
    {\cal W}_2^- + {\cal W}_3 + {\cal W}_5$.
    \item[{\rm (ii)}] If the $\Lie{G}_2$-structure on $\overline{M}$
    is of type ${\cal X}_2 + {\cal X}_4$ and $p_{\varphi}$ is tangent
     to $M$, then the $\SU(3)$-structure on $M$ is of type ${\cal W}^+_2
     + {\cal W}_2^- + {\cal W}_3 + {\cal W}_4 + {\cal W}_5$
     if and only if $M$ is a minimal variety.
       \item[{\rm (iii)}] If the $\Lie{G}_2$-structure on $\overline{M}$
    is of type ${\cal X}_2 + {\cal X}_4$, then the $\SU(3)$-structure on $M$ is of type ${\cal W}^+_1
    + {\cal W}^-_1 + {\cal W}^+_2
     + {\cal W}_2^- + {\cal W}_4 + {\cal W}_5$
     if and only if $IB=B$.
\end{enumerate}
\end{theorem}

\section{Examples} \label{sec:examples}

\noindent 1. {\it $\SU(3)$-structures on $\Lie{S}^6$}.-
 Let us consider ${\Bbb R}^7$ as identified  with   $\Lie{Im} {\Bbb O}$,
 the pure imaginary Cayley numbers. It is well known that ${\Bbb R}^7$,
considered in this way, is equipped with a parallel (${\cal P}$)
$\Lie{G}_2$-structure defined by means of the product of Cayley
numbers. Since the six-dimensional sphere $\Lie{S}^6$ is totally
umbilic in ${\Bbb R}^7$, taking Theorem~\ref{xunogeneral} into
account, $\Lie{S}^6$ has  $\SU(3)$-structures of types ${\cal
W}_1^+ + {\cal W}_1^- + {\cal W}_5$ and  ${\cal W}_1^+ + {\cal
W}_1^-$ . Likewise, taking Theorem~\ref{xunopimedio}, $\Lie{S}^6$
has two $\SU(3)$-structures of type ${\cal W}_1^+$ and ${\cal
W}_1^-$, respectively.

Moreover, since $B_{\Lie{S}^6} = \langle \cdot , \cdot \rangle$,
by Proposition~\ref{rrB}, we have $r(\omega_{\Lie{S}^6}) =
\linebreak \cos \theta  \, \langle \cdot , \cdot \rangle + \sin
\theta \, \omega_{\Lie{S}^6}$. Finally, taking  Equation
\eqref{erresu} into account, we get
\begin{equation} \label{omegaeseseis}
\nabla \omega_{\Lie{S}^6} = \cos \theta \, \psi_{\Lie{S}^6+} +
\sin \theta \, \psi_{\Lie{S}^6-}.
\end{equation}
\vspace{2mm}

\noindent 2. {\it $\SU(3)$-structures on $\Lie{S}^5 \times
\Lie{S}^1$}.- Now we consider the $\SU(3)$-structure on
$\Lie{S}^6$ with $\theta = 0$ constant. On the product manifold
$\Lie{S}^6 \times \Lie{S}^1$, we define a $\Lie{G}_2$-structure by
\begin{equation} \label{primerag2}
\varphi = - \vartheta \wedge \omega_{\Lie{S}^6} +
\psi_{\Lie{S}^6+}, \quad \ast \varphi
= - \frac{1}{2}
\omega_{\Lie{S}^6} \wedge \omega_{\Lie{S}^6}  + \vartheta \wedge
\psi_{\Lie{S}^6-},
\end{equation}
where $\vartheta$ is a  Maurer-Cartan one-form on $\Lie{S}^1$.
 Since $ \nabla \omega_{\Lie{S}^6} =  \psi_{\Lie{S}^6+}$ and
$\eta_{\Lie{S}^6}=0$, taking Theorem~\ref{igual3} and
Corollary~\ref{corigual3} into account, we have
$$
d \psi_{\Lie{S}^6+} = 0, \qquad d \psi_{\Lie{S}^6-} =  2
\omega_{\Lie{S}^6}  \wedge \omega_{\Lie{S}^6}.
$$
Therefore,
\begin{gather*}
d \varphi =  3 \vartheta \wedge \psi_{\Lie{S}^6+} =  3 \vartheta
\wedge \varphi,\\
d \ast \varphi = - 2 \vartheta \wedge \omega_{\Lie{S}^6}  \wedge
\omega_{\Lie{S}^6} =  4 \vartheta \wedge \ast \varphi.
\end{gather*}
Hence,  $\Lie{S}^6 \times \Lie{S}^1$ is equipped with a
$\Lie{G}_2$-structure of type ${\cal X}_4$ with $p d^* \varphi_1 =
- 12 \vartheta$ (see \cite{Cabrera:G2struct}). This fact was
pointed out in \cite{Cabrera:hipspin} and, moreover, we have
$\overline{r}(\varphi) =  \varphi ( \vartheta , \cdot, \cdot)$
(see Table \ref{tab:g2types}).

The manifold $\Lie{S}^5 \times \Lie{S}^1$ is contained in
$\Lie{S}^6 \times \Lie{S}^1$ as a totally geodesic orientable
hypersurface. Thus, if $n$ is the unit normal vector field on
$\Lie{S}^5 \times \Lie{S}^1$, we are in the conditions of
Corollary~\ref{xcuatrocor} with $  Id^* \omega = 12  \eta + 4 Id
\theta  = - 4 \iota^* \vartheta  \neq 0$. Therefore, the manifold
$\Lie{S}^5 \times \Lie{S}^1$ has  $\SU(3)$-structures of type
$$
\left( {\cal W}_4 + {\cal W}_5 \right) - \left( {\cal W}_4 \cup
{\cal W}_5 \right).
$$
\begin{remark}
{\rm It is well known that $\Lie{S}^6$  and $\Lie{S}^5 \times
\Lie{S}^1$ have almost Hermitian structures ($\Un(3)$-structures)
of type ${\cal W}_1$ and ${\cal W}_4$, respectively. Here we give
further information about $\Lie{S}^6$  and $\Lie{S}^5 \times
\Lie{S}^1$ as special almost Hermitian manifolds. }
\end{remark}
\vspace{2mm}

\noindent 3.{\it Hypersurfaces of the manifold $H(1,2) / \Gamma
\times {\Bbb T}^{\;2}$}.- Let $H(1,2)$ the generalized Heisenberg
group, i.e., the connected, simply connected and nilpotent Lie
group consisting of matrices: $$ a = \left( \begin{array}{ccc}
\mbox{\rm I}_2 & X & Z \\
0 & 1 &  y \\
0 & 0 & 1
\end{array} \right),
$$
where $X$ and $Z$ are $2 \times 1$ matrices of real numbers and
$y$ is a  real number. If we write $X^t = (x_1,x_2)$ and $Z^t =
(z_1,z_2)$, then a coordinate system on $H(1,2)$ is given by:
$$
x_1 (a) = x_1, \;\;\;x_2 (a) = x_2, \;\;\; y (a) = y, \;\;\;
z_1(a) = z_1, \;\;\; z_2(a) = z_2; \;\;\;
$$
and a basis of left invariant one-forms is:
$$
\{ dx_1 ,\;\; dx_2 , \;\;dy ,\;\; dz_1-x_1 dy , \;\; dz_2 - x_2 dy
\}.
$$

Let $\Gamma $ be the discrete subgroup of $H(1,2)$ consisting of
matrices with integer entries. We consider the quotient space
$H(1,2) / \Gamma $. Since the one-forms $ \{ dx_1 , dx_2 , dy,
dz_1-x_1dy , dz_2 - x_2 dy \} $ are left invariant, they descend
to one-forms $\{ \beta_1 , \beta_2, \lambda , \gamma_1 , \gamma_2
\}$ on $H(1,2) / \Gamma $ such that:
$$
\begin{array}{c}
d\beta_1 = d\beta_2 = d\lambda = 0, \qquad d\gamma_1 = \lambda
\wedge \beta_1, \qquad d\gamma_2 = \lambda \wedge \beta_2.
\end{array}
$$

Let us consider the product manifold $\overline{M}= H(1,2)/ \Gamma
\times {\Bbb T}^{\;2}$, where ${\Bbb T}^{\;2}$ is a
two-dimensional torus. If $\eta_1 , \eta_2$ denote a basis of
closed one-forms on ${\Bbb T}^{\;2}$, we consider on
$\overline{M}$ the following basis for one-forms
$$
\begin{array}{ccccccc}
e_0 = \beta_1, &  e_1 = \gamma_2 , & e_2 = \eta_2, & e_3 = \eta_1,
& e_4 = \beta_2, & e_5 = \lambda , & e_6 = \gamma_1.
\end{array}
$$

 In \cite{Fernandez:G2cal} it is defined  the $\Lie{G}_2$-structure  such that
$\left\{ e_0 , e_1 , \ldots , e_6 \right\} $ is a Cayley coframe
for one-forms and found that such a $\Lie{G}_2$-structure is of
type ${\cal X}_2$, i.e., $d \varphi = 0$. In fact, it was the
first known example of compact calibrated manifold with
$\Lie{G}_2$-structure. Moreover, one can compute
$$
d^* \varphi = - \ast d \ast \varphi = e_0 \wedge e_1 - e_4 \wedge
e_6.
$$
Now, using Equation \eqref{rxy}, we obtain the bilinear form
$$
2 \overline{r}(\varphi) = e_0 \wedge e_1 - e_4 \wedge e_6.
$$

Let us recall that we have $de_1 = - e_4 \wedge e_5$, $d e_6 = -
e_0 \wedge e_5$ and $de_i = 0$, for $i=0,2,3,4,5$. From $e_i ( [
e_j , e_k ] ) = - d e_i (e_j,e_k)$, where $i,j,k \in {\Bbb
Z}_{\;7}$ and $[,]$ denotes the Lie bracket, we get $[e_0,e_5] =
e_6$, $[e_4 ,e_5] =  e_1$ and $[e_i,e_j]=0$ for $(i,j) \in {\Bbb
Z}_{\,7} \times {\Bbb Z}_{\,7} -  \{ (0,5), (5,0), (4,5) , (5,4)
\}$. Now, using Koszul Formula,  we get $\overline{\nabla}_{e_0}
e_5 = \frac{1}{2} e_6$, $\overline{\nabla}_{e_0} e_6 = -
\frac{1}{2} e_5$,
 $\overline{\nabla}_{e_1} e_4 = - \frac{1}{2} e_5$,
 $\overline{\nabla}_{e_1} e_5 =  \frac{1}{2} e_4$
 $\overline{\nabla}_{e_4} e_1 = - \frac{1}{2} e_5$,
 $\overline{\nabla}_{e_4} e_5 = - \frac{1}{2} e_1$,
 $\overline{\nabla}_{e_5} e_0 = - \frac{1}{2} e_6$,
 $\overline{\nabla}_{e_5} e_1 =  \frac{1}{2} e_4$,
 $\overline{\nabla}_{e_5} e_4 = - \frac{1}{2} e_1$,
 $\overline{\nabla}_{e_5} e_6 =  \frac{1}{2} e_0$,
 $\overline{\nabla}_{e_6} e_0 = - \frac{1}{2} e_5$,
 $\overline{\nabla}_{e_6} e_5 =  \frac{1}{2} e_0$,
and $\overline{\nabla}_{e_i} e_j = 0,$ for all  $(i,j) \in {\Bbb
Z}_{\;7} \times {\Bbb Z}_{\;7} - \{ (0,5), (0,6), (1,4), (1,5),
(4,1), (4,5), (5,0), (5,1), (5,4), (5,6), (6,0), (6,5) \}$.
\vspace{3mm}

\begin{enumerate}
\item[(i)] Let us consider a hypersurface $M_1$ which is a maximal
integral submanifold of the integrable distribution defined by the
one-form $n=e_3$. Such hypersurface $M_1$ is diffeomorphic with $M
\times {\Bbb T}^{\;1}$.
 Since  $M_1$ is totally geodesic and  $\overline{r}(\varphi)(\iota_*
 \cdot , n)=0$, making use of Theorem~\ref{xdosgeneral}, $M_1$ has
 induced an $\SU(3)$-structure of type ${\cal W}_2^+ + {\cal W}_2^- + {\cal
 W}_5$. If we take $\theta = \pi / 4$ constant, the $\SU(3)$-structure
 is of type ${\cal W}_2^+ + {\cal W}_2^-$.  Finally, making use of
 Theorem~\ref{xdospimedio}, we obtain  two  $\SU(3)$-structures on $M_1$
 of type ${\cal W}_2^+$ and ${\cal W}_2^-$, respectively.

\item[(ii)]  Let us consider a hypersurface $M_2$ which is a
maximal integral submanifold of the integrable distribution
defined by the one-form $n=e_0$. The second fundamental form $M_2$
is  given by $ 2B = e_5 \otimes e_6 + e_6 \otimes e_5$ and $2
\iota^* \overline{r}(\varphi) = - e_4 \wedge e_6$. Thus, $M_2$ is
a minimal variety and it is satisfied $(I_{(1)} - I_{(2)}) \iota^*
\overline{r}(\varphi) + (1+I)B=0$. Therefore, by making use of
Theorem~\ref{xdosgeneral}, the $\SU(3)$-structure on $M_2$ is of
type ${\cal W}_3 + {\cal W}_4 + {\cal W}_5$. Note that $IB \neq B$
and $-2 \overline{r}(\varphi)(\iota_* \cdot , n ) = e_1 \neq 0$.

 \item[(iii)]Let us consider a hypersurface $M_3$ which is a
maximal integral submanifold of the integrable distribution
defined by the one-form $n=e_5$. The second fundamental form of
$M_3$ is  given by $ B = - e_0 \vee e_6 - e_1 \vee e_4$, where $2
a \vee b= a \otimes b + b \otimes a$. Moreover, $2 \iota^*
\overline{r}(\varphi) = e_0 \wedge e_1 - e_4 \wedge e_6$. Thus,
$M_3$ is a minimal variety, $IB = B $ and
$\overline{r}(\varphi)(\iota^* \cdot , n)= 0$. Then, using
Theorem~\ref{xdosgeneral}, we have that $M_3$ is of type
K{\"a}hler (${\cal W}_5$). In particular, when $\theta$ is
constant, the intrinsic $\SU(3)$-torsion of $M_3$ vanishes, i.e.,
$M_3$ is $\su(3)$-K{\"a}hler.
\end{enumerate}
\vspace{2mm}

\noindent 4.{\it Hypersurfaces of a family of manifolds with
$\Lie{G}_2$-structure of type ${\cal X}_3$}.- Let us consider the
manifolds $M(k)$ described in \cite{CFG:example} as follows. For a
fixed $k \in \mathbb R , \; k\neq 0$, let $G(k)$ be the
three-dimensional connected and solvable (non-nilpotent) Lie group
consisting  of the matrices
$$
{\bf a} = \left( \begin{array}{cccc}
e^{k z} & 0 & 0 & x \\
0 & e^{-kz} & 0 & y \\
0 & 0 & 1 & z \\
0 & 0 & 0 & 1
\end{array} \right),
$$
where $x,y,z \in \mathbb R$. Then, a global coordinate system $\{
x , y , z \}$  for $G(k)$ is given by $ x({\bf a})=x, \; y({\bf
a}) = y, \; z({\bf a}) = z . $ A straightforward computation
proves that a basis of right invariant one-forms on $G(k)$ is $\{
 dx - k x dz , dy + k y dz ,  dz \}$.

The Lie group $G(k)$ can be also described as the semidirect
product $\mathbb R \times_{\phi} \mathbb R^{2}$, where $ \phi \, :
\, \mathbb R \, \to \, \Lie{Aut}(\mathbb R^{2}) $ is the
representation defined by
$$
\phi (t) = \left( \begin{array}{cc}
e^{k z} & 0 \\
0 & e^{-k z}
\end{array} \right) ,\;\;\; z \in \mathbb R.
$$
Therefore $G(k) $ possesses a discrete subgroup $\Gamma (K)$ such
that the quo\-tient ma\-ni\-fold $M(k) = G(k) / \Gamma (k)$ is
compact. Moreover, the one-forms $dx - k x dz, \; dy + k y dz,\;
dz$ descend to $M(k)$. Let us denote by $ \alpha, \beta , \gamma$,
respectively, the induced one-forms on $M (k)$. Then  we have $
d\alpha = - k \alpha \wedge \gamma , \;\;  d\beta = k \beta \wedge
\gamma , \;\;  d \gamma = 0. $

Let us consider the product manifold $M (k)
 \times {\Bbb T}^{\;4}$ where ${\Bbb T}^{\;4}$ is
a four-dimensional torus. Let $\eta_{1} , \eta_{2} , \eta_{3} ,
\eta_{4}$ be a basis of closed one-forms on ${\Bbb T}^{\;4}$. Then
in $M (k) \times {\Bbb T}^{\;4}$ we consider the following basis
for  one-forms
$$
\begin{array}{ccccccc}
   e_{0} = \alpha, & e_{1} =
\beta, & e_{2} = \eta_{2}, & e_{3} = \gamma , & e_{4} = \eta_{3},
& e_{5} = \eta_{4}, & e_{6} = \eta_{5}.
\end{array}
$$
Therefore,  $ d e_{0} =  - k  e_{0} \wedge e_{3}$, $d e_{1} = k
e_{1} \wedge e_{3}$ and $d e_i=0$, for all $i \in {\Bbb Z}_{\;7} -
\{ 0,1 \}$.

We consider on $M (k) \times {\Bbb T}^{\;4}$ the
$\Lie{G}_2$-structure determined by the four  form
$
 \varphi  =  \sum_{i \in {\Bbb  Z}_{\; 7}}  e_{i}  \wedge e_{i+1}  \wedge
 e_{i+3}.
$
Since
\begin{eqnarray} \label{ejemplo2}
 d \varphi &  = & \;  k \,  e_{3} \wedge e_{4} \wedge e_{5}
\wedge e_0 - k \,  e_{3} \wedge
e_{5}  \wedge  e_{6} \wedge e_1 \\
 & & + k  e_{3}  \wedge e_{6}  \wedge  e_{0} \wedge e_2   - k \,
e_{3}  \wedge e_{1}  \wedge e_{2} \wedge e_4, \nonumber
\end{eqnarray}
we have $d \varphi \wedge \varphi =0$. Moreover, it is immediate
to check that $d \ast \varphi = 0$. Therefore,  the ${\rm
G}_2$-structure is of type ${\cal X}_3$. Moreover, with the
indicated metric, these manifolds do not admit any
$\Lie{G}_2$-structure of type ${\cal P}$ (see
\cite{Cabrera:hipspin}).

Now, from $d^* \varphi =0 $ and the expression for $d \varphi$
given by \eqref{ejemplo2}, taking Equation \eqref{rxy} into
account, we obtain
\begin{equation} \label{ejem2t}
r(\overline{\varphi}) = - 2 k  e_0 \vee e_1.
\end{equation}

From $e_i ( [ e_j , e_k ] ) = - d e_i (e_j,e_k)$, for $i,j,k \in
{\Bbb Z}_{\;7}$, we get $[e_0,e_3] = k e_0$, $[e_1 ,e_3] = - k
e_1$ and $[e_i,e_j]=0$ for $(i,j) \in {\Bbb Z}_{\,7} \times {\Bbb
Z}_{\,7} -  \{ (0,3), (3,0), (1,3) , (3,1) \}$. Now, using Koszul
Formula,  we get $\overline{\nabla}_{e_0} e_3 = k e_0$,
 $\overline{\nabla}_{e_1} e_3 = - k e_1$,
 $\overline{\nabla}_{e_0} e_0 = - k e_3$,
 $\overline{\nabla}_{e_1} e_1 = k e_3$
and $\overline{\nabla}_{e_i} e_j = 0,$ for all  $(i,j) \in {\Bbb
Z}_{\;7} \times {\Bbb Z}_{\;7} - \{ (0,3), (1,3), \linebreak
(0,0), (1,1) \}$.
\begin{enumerate}
\item[(i)] Let us consider a hypersurface $N_1$ which is a maximal
integral submanifold of the integrable distribution defined by the
one-form $n=e_2$. Such a hypersurface $N_1$ is diffeomorphic with
$M(k) \times {\Bbb T}^{\;3}$.
 Since  $N_1$ is totally geodesic,
  $r(\overline{\varphi}) (  \iota_* \cdot , n ) =0$, $r(\overline{\varphi}) ( n , n) =0$,  and
 $\iota^*r(\overline{\varphi}) = - 2 k  e_0 \vee e_1$, taking
 Theorem~\ref{xtresgeneral} into account,  $N_1$
 has  $\SU(3)$-structure of types ${\cal W}_2^+ + {\cal W}_2^-
 + {\cal W}_3 + {\cal W}_5$. Taking $\theta$ constant, the $\SU(3)$-structure on $N_1$ is of type
 ${\cal W}_2^+ + {\cal W}_2^- + {\cal W}_3 $.  In particular, when
 $\theta=0(\theta=\pi/2)$ constant, making use of
 Theorem~\ref{xtrespimedio}, we have that the $\SU(3)$-structure
 is of type ${\cal W}_2^+ + {\cal W}_3$ (${\cal W}_2^- + {\cal
 W}_3$).

\item[(ii)] Let us consider a hypersurface $N_2$ which is a
maximal integral submanifold of the integrable distribution
defined by the one-form $n=e_3$. The second fundamental tensor of
$N_2$ is given by
 $
B = k e_0 \otimes e_0 - k e_1 \otimes e_1$.  Thus, $M_2$ is a
minimal variety, $r(\overline{\varphi}) ( \iota_* \cdot , n ) =0$,
$r(\overline{\varphi}) ( n , n ) =0$, $IB=-B$ and  $I
\iota^*r(\overline{\varphi}) = - \iota^*r(\overline{\varphi})$.
Then, using Theorem~\ref{xtresgeneral},  $N_1$
 has  $\SU(3)$-structures of type ${\cal W}_3 + {\cal W}_5$.
In particular, when $\theta$ is constant, the $\SU(3)$-structure
on $N_2$ is of
 type ${\cal W}_3$.
\item[(iii)]  Let us consider a hypersurface $N_3$ which is a
maximal integral submanifold of the integrable distribution
defined by the one-form $n=e_0$. It is immediate that $N_3$ is
totally geodesic, $r(\overline{\varphi}) (  \iota_* \cdot , n ) =-
ke_1$, $\iota^*r(\overline{\varphi})=0$ and $r(\overline{\varphi})
( n , n ) =0$. Then, using Theorem~\ref{xtresgeneral}, the
$\SU(3)$-structure on $N_3$ is of type ${\cal W}_5$. If we
consider $\theta$ constant, then $N_3$ is
$\lie{su}(3)$-K{\"a}hler.
\end{enumerate}

\begin{footnotesize}
 \setlength{\parindent}{0pt} \setlength{\parskip}{10pt}

  \textit{Department of Fundamental Mathematics,
  University of La Laguna, 38200 La Laguna, Tenerife, Spain}.  E-mail:
  {\tt fmartin@ull.es}

\end{footnotesize}

\end{document}